%% file: arxiv_lp_Nov2026.tex
\documentclass[11pt]{article}

\usepackage{geometry}

\usepackage{url}
\usepackage{amssymb,enumerate}
\usepackage{amsmath,amsfonts}
\usepackage{amsthm}
\usepackage{latexsym}
\usepackage{mathrsfs}
\usepackage{color}
\usepackage{xcolor}
\usepackage{microtype}

\usepackage{algorithm}
\usepackage{algorithmicx}
\usepackage{algpseudocode}
\usepackage{booktabs}
\usepackage{enumitem}
\usepackage[most]{tcolorbox}
\usepackage[normalem]{ulem}
\usepackage{multirow}
\allowdisplaybreaks[4]
\usepackage[numbers,square,sort]{natbib}
\usepackage[
  colorlinks=true,
  linkcolor=teal,
  citecolor=teal,
  urlcolor=teal
]{hyperref}

\input{math_commands.tex}

\theoremstyle{plain}
\newtheorem{theorem}{Theorem}

\newtheorem{lemma}{Lemma}

\newtheorem{assumption}{Assumption}

\newtheorem{remark}{Remark}

\usepackage{graphicx, color}
\usepackage{algorithm, algpseudocode} % use algorithm and algorithmicx for typesetting algorithms
\usepackage{mathrsfs} % for \mathscr command

\usepackage{lipsum}

\makeatletter
\newenvironment{breakablealgorithm}
  {% \begin{breakablealgorithm}
   \begin{center}
     \refstepcounter{algorithm}
     \hrule height.8pt depth0pt \kern2pt
     \renewcommand{\caption}[2][\relax]{
       {\raggedright\textbf{\ALG@name~\thealgorithm} ##2\par}%
       \ifx\relax##1\relax
         \addcontentsline{loa}{algorithm}{\protect\numberline{\thealgorithm}##2}%
       \else
         \addcontentsline{loa}{algorithm}{\protect\numberline{\thealgorithm}##1}%
       \fi
       \kern2pt\hrule\kern2pt
     }
  }{% \end{breakablealgorithm}
     \kern2pt\hrule\relax
   \end{center}
  }
\makeatother
\newcounter{calibration}

\title{A simple adaptive accelerated optimistic dual averaging method for unconstrained convex optimization in $p$-norm geometry}
\author{
Bowen Li\thanks{Department of Mathematics, Link\"oping University, Sweden
(Email: \texttt{bowen.a.li@liu.se},
\texttt{danqing.zhou@liu.se}, and
\texttt{chuan.he@liu.se}). All authors contributed equally to this work.}
\and
Danqing Zhou\footnotemark[1]
\and
Chuan He\footnotemark[1]
}

\date{
}

\begin{document}
\maketitle
    
\begin{abstract}
Efficient preconditioning plays an important role in large-scale first-order optimization. In this paper, we propose an adaptive accelerated optimistic dual averaging ($\text{A}^2$ODA) method for unconstrained convex optimization that adopts low-cost preconditioning based on $p$-norm geometry for $p\in[2,\infty]$. The method is line-search-free and employs a new stepsize selection technique based on an aggregated affine lower model. Its update structure is simple and unified across $p\in[2,\infty]$, with the case $p>2$ requiring no more elaborate updates than the Euclidean case $p=2$. Under a generalized local smoothness condition, we derive a
convergence bound that separates the effects of constant and
gap-dependent curvature, with the latter potentially yielding
faster-decaying terms. In the special case of local smoothness,
this bound recovers a nearly optimal convergence rate.

\end{abstract}


\smallskip


\section{Introduction}
In this paper, we consider the unconstrained convex optimization problem
\begin{equation}\label{eq:problem}
    \min_{x\in\mathbb R^n} f(x),
\end{equation}
where $f:\mathbb R^n\to\mathbb R$ is a continuously differentiable convex function. We assume that problem \eqref{eq:problem} has at least one optimal solution. Since the seminal work on accelerated gradient methods~\cite{nesterov1983method}, the development of fast first-order methods has remained a central focus in the study of smooth convex optimization. In practice, the performance of first-order methods is influenced by many factors, including the acceleration mechanism, the parameter selection, and the underlying 
geometry. The goal of this paper is to develop an accelerated first-order method in $p$-norm geometry that adapts its parameters to the
optimization trajectory and achieves a unified convergence bound
for all $p\in[2,\infty]$.

In recent years, there has been fruitful research on first-order methods for smooth convex optimization~\cite{Nesterov2004,allenzhu2017linear,nesterov2015universal}. In particular, many recent efforts have focused on parameter selection in first-order methods, which is often guided by estimates of the Lipschitz constant of the objective gradient. Instead of fixing the parameters based on a global Lipschitz constant, adapting them to local curvature typically leads to much faster practical progress. In this direction, there has recently been growing interest in first-order methods that adapt their parameters based on local curvature without requiring costly line searches. The use of local curvature to select stepsizes has a long history, as exemplified by the well-known Barzilai-Borwein (BB) stepsize~\cite{BB88}. However, directly applying the BB stepsize without line searches does not generally yield convergence guarantees for smooth convex problems. Recently, several studies have proposed novel stepsize schemes for line-search-free gradient methods and established provable convergence guarantees for solving problem \eqref{eq:problem}; see, e.g., \cite{MalitskyMishchenko2020,LTP2024,ZMY2026}. In addition, notable adaptive accelerated gradient methods, including AC-FGM~\cite{li2026simple} and AdaNAG~\cite{suh2025adaptive}, achieve optimal convergence rates for smooth convex optimization without requiring line searches.

Beyond parameter selection, the choice of geometry is another important aspect of first-order method design. For many constrained problems, replacing the Euclidean squared distance with a suitable Bregman divergence in first-order method design allows the resulting updates to reflect the intrinsic geometry of the feasible region. Well-known examples include the Kullback-Leibler divergence on the probability simplex~\cite{beck2003mirror} and the Itakura-Saito divergence on the positive orthant~\cite{bauschke1997legendre}. These examples can be unified under the general mirror descent framework \cite{nemirovski1983problem,beck2003mirror}, whose update takes the form:
\begin{equation*}
x_{k+1}\in\operatorname*{arg\,min}_{x\in X}
\left\{\langle\nabla f(x_k),x-x_k\rangle+\gamma_kD_h(x,x_k)\right\},
\end{equation*}
where $\gamma_k>0$, $X$ denotes the feasible region, and
$D_h(x,y):=h(x)-h(y)-\langle\nabla h(y),x-y\rangle$
is the Bregman divergence generated by $h$. When \(h(x)=\|x\|_2^2/2\), the update reduces to projected gradient descent. More generally, suitable choices of \(h\) allow mirror descent to exploit non-Euclidean geometry. Under the relative smoothness condition
\begin{equation*}%\label{rela-smooth}
D_f(x,y)\le L D_h(x,y)\qquad\forall x,y\in\Omega
\end{equation*}
for some $L>0$, convergence bounds of mirror descent methods can be established in terms of \(D_h\), where \(\Omega\) is a suitable convex set containing the iterates and relevant comparison points~\cite{beck2003mirror,lu2018relative}. Mirror descent and relative smoothness extend gradient descent and Lipschitz smoothness in the Euclidean norm, respectively, to non-Euclidean settings. Nonetheless, extending accelerated gradient methods to non-Euclidean mirror descent is nontrivial. There have been many recent studies on acceleration in mirror descent, including both specific accelerated variants~\cite{krichene2015accelerated, hanzely2021accelerated, MST2026} and lower bounds on achievable convergence rates~\cite{dragomir2022optimal}. In addition, adaptive mirror descent methods have also been developed~\cite{ou2026localrelative}.

Another line of research explores first-order methods that directly exploit norm geometry to shape search directions. In particular, replacing the squared Euclidean distance with the squared distance induced by a
general norm yields the dual preconditioned gradient update
\cite{kelner2014almost}:
\begin{align*}
x_{k+1}
= x_k-\frac{1}{\gamma_k}J_*(\nabla f(x_k))
\in\operatorname*{arg\,min}_{x\in\R^n}
\left\{
\langle \nabla f(x_k),x-x_k\rangle
+\frac{\gamma_k}{2}\|x-x_k\|^2
\right\},
\end{align*}
where $\gamma_k>0$, $\|\cdot\|$ and $\|\cdot\|_*$ denote a general norm and its dual, respectively, and $J_*(g)\in\partial(\|\cdot\|_*^2/2)(g)$ for all $g\in\R^n$. For convex problems, this preconditioned method achieves an $\mathcal{O}(T^{-1})$ convergence rate in the objective gap under a Lipschitz gradient condition defined using the general norm and its dual, together with a level-boundedness assumption \cite{kelner2014almost}, where $T$ denotes the number of iterations. This bound is analogous to that of classical gradient descent. Moreover, this preconditioned gradient descent has also been analyzed under a dual relative smoothness condition \cite{maddison2021dual}. This condition serves as a dual counterpart to the primal relative smoothness condition used to analyze mirror descent \cite{lu2018relative}. In addition, accelerated variants of preconditioned gradient descent in general norm geometries have been developed in \cite{diakonikolas2024complementary,daspremont2018affine}. In particular, when the primal and dual norms are the $p$- and $q$-norms, respectively, with $1/p+1/q=1$ and $2\le p<\infty$, these methods achieve an $\mathcal{O}(T^{-(1+2/p)})$ convergence
rate in the objective gap, matching the lower bounds up to logarithmic factors
\cite{diakonikolas2020lower}. From a technical standpoint, these methods rely on the uniform convexity of
the proximal regularizer to control linearization errors in their accelerated updates \cite{diakonikolas2024complementary,daspremont2018affine}. For the case $p>2$, the resulting $p$th-power proximal terms interact with the quadratic terms arising from Lipschitz smoothness, making parameter selection and the convergence analysis more involved than in the case $p=2$ (see Remark \ref{rem:relaxed-dual-expansion}). It is also noteworthy that $p$-norm geometry has attracted significant recent interest
in the design of deep learning optimizers, particularly in the
endpoint case $p=\infty$, which leads to sign-based update directions; see, e.g.,
\cite{jordan2024muon,bernstein2018signsgd,dong2026muonp,ma2026schattor,pethick2026schatten}. These theoretical advances and growing practical interest motivate the following natural question:
\begin{center}
\emph{Can we design a line-search-free accelerated first-order method for
smooth convex optimization that exploits $p$-norm geometry for
preconditioning throughout the range $p\in[2,\infty]$?} 
\end{center}

This paper provides an affirmative answer to this question.
Specifically, our proposed method employs an acceleration mechanism
that combines optimistic dual averaging~\cite{rakhlin2013online}
in $p$-norm geometry with weighted iterate averaging~\cite{joulani2020averaging}.
To adapt the regularization parameters, we aggregate past linearizations
into an affine model and use this model to derive a new parameter
selection rule that directly controls a Lyapunov potential. This construction differs fundamentally
from the parameter selection rules used by existing line-search-free
accelerated gradient methods in the case $p=2$~\cite{li2026simple,suh2025adaptive}.

Moreover, our analysis relies on the relaxed dual-potential expansion
inequality in Lemma~\ref{lem:expan-vq}, which avoids the need for
a uniform convexity inequality for the proximal regularizer.
Consequently, our method directly uses the squared $p$-norm as
the proximal regularizer for $p>2$, in contrast to existing accelerated
methods that rely on $p$-uniformly convex
regularizers~\cite{diakonikolas2024complementary,daspremont2018affine}
(see Remark~\ref{rem:relaxed-dual-expansion}). On the theoretical side, we establish an accelerated convergence rate for the proposed method under a generalized local smoothness condition (see Assumption~\ref{asp:cvx}(b)). This condition incorporates gap-dependent curvature into the local smoothness condition and allows us to distinguish the contributions
of the constant and gap-dependent curvature terms to the convergence bound. Our main contributions are summarized below.

\begin{itemize}[leftmargin=1.2em]
    \item We develop an adaptive accelerated optimistic dual averaging ($\text{A}^2$ODA) method under $p$-norm geometry for solving problem~\eqref{eq:problem}. To the best of our knowledge, this is the first line-search-free first-order method that combines adaptive acceleration with squared $p$-norm regularization in a unified framework covering all $p\in[2,\infty]$.

    \item Under the generalized local smoothness condition
(Assumption~\ref{asp:cvx}(b)), we establish a unified convergence
bound for Algorithm~\ref{alg:pODA} across all $p\in[2,\infty]$
that explicitly separates the contributions of constant and
gap-dependent curvature
(see Theorem~\ref{thm:convex-convergence-rate}).
In the special case of local smoothness, this bound recovers
a nearly optimal convergence rate
(see Section~\ref{subsec:complete-alg}).
\end{itemize}

The rest of the paper is organized as follows. In Section~\ref{sec:preliminaries}, we introduce the notation and preliminary results. In Section~\ref{sec:algorithm}, we present the proposed algorithm
and state its convergence guarantees.  The proof of our main convergence result
is presented in Section~\ref{sec:proof}. Additional lemmas and proofs are deferred to the appendix.

\section{Notation and preliminaries}\label{sec:preliminaries}

%\paragraph{Notation.} 

Throughout this paper, we use \(\mathbb{R}^n\) and $\R_+$ to denote the \(n\)-dimensional Euclidean space and the set of nonnegative real numbers, respectively, and \(\langle\cdot,\cdot\rangle\) to denote the standard inner product. %$e$ denotes Euler's number 
We adopt the convention $0^0=1$ and $1/\infty=0$. For a positive integer \(N\), we define $[N]:=\{1,2,\ldots,N\}$. For any \(p\in[1,\infty)\) and \(z\in\mathbb{R}^n\), we define $\|z\|_p:=(\sum_{i=1}^n |z_i|^p)^{1/p}$ and $\|z\|_\infty:=\max_{1\le i\le n}|z_i|$.  For any $A\in\R^{n\times n}$ and any $p,q\in[1,\infty]$, the induced $p$-to-$q$ operator norm is defined by $\|A\|_{p\to q} := \sup_{\|x\|_p\le1}\|Ax\|_q$. For any convex function $h:\R^n\to\R$, its Fenchel conjugate is defined by $h^*(u):=\sup_{x\in\R^n}\{\langle x,u\rangle - h(x)\}$ for every $u\in\R^n$. Let $\partial h(x)$ denote the convex subdifferential of $h$ at $x$. The Fenchel-Young equality holds if and only if $u\in\partial h(x)$:
\begin{align}\label{eq:FY}
    h(x)+h^*(u)=\langle x,u\rangle.
\end{align}
For $x,y\in\mathbb{R}^n$ and $g_y\in\partial h(y)$, define the generalized Bregman divergence by $D_h^{g_y}(x,y):=h(x)-h(y)-\langle g_y,x-y\rangle$. When $h$ is differentiable at $y$, we omit the superscript and denote the Bregman divergence by $D_h(x,y)$. In addition, we use $\mathcal{O}(\cdot)$ and $\Theta(\cdot)$ to denote
asymptotic upper and two-sided bounds, respectively.

%In addition, we use $\mathcal{O}(\cdot)$ to denote the standard big-O notation.
For conjugate exponents $p\in[2,\infty]$ and $q\in[1,2]$ satisfying $1/p+1/q=1$,
% \footnote{We adopt the convention $1/\infty=0$.}, 
define the primal and dual regularizers by
\begin{align}\label{def:V}
    d_p(x) :=\frac12\|x\|_p^2, \qquad V_q(u) :=d_p^*(u) =\frac12\|u\|_q^2 \qquad \forall x,u\in\mathbb{R}^n.
\end{align}
For any $q\in[1,2]$ and $u\in\mathbb{R}^n$, define
\begin{align}\label{def:J}
    J_q(u) &:= \begin{cases}
        \nabla V_q(u), & q\in(1,2],\\
        \|u\|_1\operatorname{sgn}_0(u), & q=1,
    \end{cases}
    &
    (\operatorname{sgn}_0(u))_i &:= \begin{cases}
        1,  & u_i>0,\\
        0,  & u_i=0,\\
        -1, & u_i<0,
    \end{cases}
    \qquad i=1,\ldots,n.
\end{align}
Then, for any pair of conjugate exponents $p\in[2,\infty]$ and $q\in[1,2]$, one can verify that
\begin{align}\label{ppty:J}
    J_q(u)\in\partial V_q(u),
    \qquad
    \|J_q(u)\|_p=\|u\|_q,
    \qquad
    \langle J_q(u),u\rangle=\|u\|_q^2.
\end{align}
For any $q\in[1,2]$, we also introduce an exponent that plays an important role in our analysis:
\begin{align}\label{def:rq}
    r_q := \begin{cases}
        \frac{q}{2-q},  & q\in[1,2),\\
        \infty,  & q=2.
    \end{cases} 
\end{align}

The next lemma establishes an expansion inequality for the
squared $\ell_q$-norm, including the nonsmooth case $q=1$.
It provides a counterpart of the well-known Euclidean identity
\begin{align}\label{eq:E-idt}
    V_2(u+v) = V_2(u)+\langle J_2(u),v\rangle+V_2(v).
\end{align}

\begin{lemma}[Relaxed dual-potential expansion]\label{lem:expan-vq}
Let $q\in[1,2)$, and let $V_q$ and $J_q$ be defined in
\eqref{def:V} and \eqref{def:J}, respectively. Then,
\begin{align}\label{ineq:V-expan}
c\,V_q(u+v) \le V_q(u) + \langle J_q(u),v\rangle + \bigg(1+\frac{2^{2/r_q}c}{(1-c^{r_q})^{1/{r_q}}}\bigg)\,V_q(v)
\end{align}
holds for all $u,v\in\R^n, c\in(0,1)$, where $r_q$ is defined in \eqref{def:rq}.
\end{lemma}

\begin{remark}
\label{rem:relaxed-dual-expansion}
Existing accelerated methods in $p$-norm geometry typically rely on $p$-uniformly convex proximal regularizers for $p>2$~\cite{diakonikolas2024complementary,daspremont2018affine}. The resulting $p$th-power proximal terms cannot be directly balanced against the quadratic errors arising from Lipschitz smoothness, requiring additional approximation arguments and more involved parameter choices.

In contrast, our method directly uses the squared $p$-norm as the proximal regularizer, without requiring a uniform convexity lower bound. The key technical ingredient in our analysis is the relaxed dual-potential expansion in Lemma~\ref{lem:expan-vq}, which controls perturbations through a quadratic term in the dual norm. Combined with our adaptive parameter-selection rule, this inequality provides an alternative approach to acceleration in $p$-norm geometry and simplifies the convergence analysis.
\end{remark}

We make the following assumptions throughout the paper.

\begin{assumption}\label{asp:cvx}
    \begin{enumerate}[leftmargin=*]
        \item[(a)] The function $f$ is convex and attains its minimum, with optimal value denoted by $f_*:=\min_{x\in\mathbb{R}^n}f(x)$.

        \item[(b)] Let $p\in[2,\infty]$ and $q\in[1,2]$ be conjugate exponents with $1/p+1/q=1$, and let $\nu\in[0,1)$. For every compact convex set $\Omega\subset\mathbb R^n$, there exist constants $L_p(\Omega),H_p(\Omega)\ge0$ such that
        \begin{align}\label{ineq:local-smoothness}
            \|\nabla f(x)-\nabla f(y)\|_q^2
            \le 2\left(L_p(\Omega)+H_p(\Omega)
            \max\{\Delta(x),\Delta(y)\}^{\nu}\right)D_f(x,y)
        \end{align}
        holds for all $x,y\in\Omega$, where $\Delta(x):=f(x)-f_*$.
    \end{enumerate}
\end{assumption}

\begin{remark}\label{rmk:local-smoothness}

%\begin{enumerate}[leftmargin=*]
    Assumption \ref{asp:cvx}(b) is a cocoercivity-type generalized local smoothness condition. This condition is related to local smoothness~\cite{lu2023locally} and gap-dependent generalized smoothness~\cite{lobanov2026gap}. Applying \eqref{ineq:local-smoothness} to both \((x,y)\) and \((y,x)\), and then using Hölder's inequality, yields the following Lipschitz-type generalized local smoothness condition:
    \begin{align}\label{ineq:general-smooth}
        \|\nabla f(x)-\nabla f(y)\|_q \le (L_{p}(\Omega) + H_p(\Omega) \max\{\Delta(x), \Delta(y)\}^\nu)\|x-y\|_p \qquad \forall x,y\in\Omega.
    \end{align}
    Conversely, for any given $\Omega$ containing $x_*$, \eqref{ineq:local-smoothness} holds with constants $L_p(\Omega),H_p(\Omega)\ge0$ provided that
    \begin{align*}
        \|\nabla f(x)-\nabla f(y)\|_q
    \le
    (
        L_p(\Omega)+2^{-\nu}H_p(\Omega)
        \max\{\Delta(x),\Delta(y)\}^{\nu}
    )\|x-y\|_p
    \qquad \forall x,y\in\bar{\Omega},
    \end{align*}
    where $\bar{\Omega}$ is a suitably enlarged compact convex set containing $\Omega$. A precise choice of $\bar{\Omega}$ and a proof of this implication are provided in Appendix~\ref{apx:cc}. Thus, the cocoercivity-type and Lipschitz-type generalized local smoothness conditions are equivalent, up to an enlargement of the domain and a change in constants.
    
    When $\nu=0$, Assumption~\ref{asp:cvx}(b) reduces to the usual local
smoothness condition. For $\nu>0$, it separates the curvature bound
into a constant part $L_p(\Omega)$ and a gap-dependent part
$H_p(\Omega)\max\{\Delta(x),\Delta(y)\}^{\nu}$, which vanishes as
the objective values approach the minimum. This decomposition allows the leading term of the convergence bound
to depend on $L_p(\Omega)$, which can be substantially smaller than
a uniform curvature bound over $\Omega$, while the gap-dependent
curvature contributes terms that may decay faster.%As will be discussed in Section~\ref{sec:applications}, for statistical applications, the constant curvature modulus \(L_p(\Omega)\) often depends only on a bound on the Hessian norm or other information of the objective function at a minimizer. This often results in improved dependence on problem constants in the leading terms of the convergence rate.
%\end{enumerate} 
    
\end{remark}

\section{An adaptive accelerated optimistic dual averaging method}\label{sec:algorithm}

In this section, we propose an adaptive accelerated optimistic dual averaging ($\text{A}^2$ODA) method in $p$-norm geometry for $p\in[2,\infty]$ and establish its convergence properties under the generalized local smoothness condition stated in Assumption~\ref{asp:cvx}(b). Specifically, we first describe the main components of the method in Sections~\ref{subsec:oda-a}-\ref{subsec:blpc} and then present the complete algorithm and its convergence results in Section~\ref{subsec:complete-alg}. 

Before proceeding, we let $\{\lambda_k\}$ be the positive weights for averaging iterates and gradients, and define the
cumulative weights by $A_k=\sum_{i=1}^k\lambda_i$ for all $k\ge1$. The $\text{A}^2$ODA method generates two sequences: the optimistic dual-averaging iterates $\{y_k\}$ and the weighted-average iterates $\{x_k\}$. For convenience, we denote the sequences of function values and gradients evaluated at $\{x_k\}$ by $\{f_k\}$ and $\{g_k\}$, respectively, where $f_k=f(x_k)$ and $g_k=\nabla f(x_k)$ for all $k\ge1$.

\subsection{Optimistic dual averaging and acceleration}\label{subsec:oda-a}

The proposed method combines optimistic dual averaging with weighted primal averaging to achieve acceleration. It mainly consists of three parts: dual averaging to accumulate past gradients, optimism to predict the next gradient contribution, and weighted primal averaging to construct the iterates at which gradients are evaluated.

Following the dual-averaging strategy~\cite{xiao2010dual}, we
maintain the weighted gradient sum
\[
    s_0=0,\qquad
    s_k=\sum_{i=1}^k\lambda_i g_i,
    \qquad \forall k\ge1.
\]
At the beginning of iteration $k$, the new gradient $g_k$ is
not yet available. We therefore use the most recent gradient
$g_{k-1}$ to predict the contribution of $g_k$ to $s_k$, forming
\[
    \tilde{s}_k=s_{k-1}+\lambda_k g_{k-1}\qquad\forall k\ge1.
\]
The term $\lambda_k g_{k-1}$ constitutes the
optimistic part: it adds a prediction of the
next weighted gradient before that gradient is observed,
as in optimistic follow-the-regularized-leader
methods~\cite{rakhlin2013online}.

The optimistic dual-averaging iterate $y_k$ is then obtained
by minimizing the sum of the predicted linear model and a
regularization term:
\[
    y_k\in\argmin_{y\in\mathbb{R}^n}
    \left\{
        \langle\tilde{s}_k,y-x_0\rangle
        +\frac{\gamma_k}{2}\|y-x_0\|_p^2
    \right\},
\]
where $\gamma_k>0$ controls the regularization strength. Our novel adaptive strategy for updating $\gamma_k$ will be introduced later.

To achieve acceleration, a key ingredient is the coupling of optimistic updates with gradient evaluations at running weighted
averages~\cite{joulani2020averaging}. Specifically, we form the weighted primal average
\[
    x_k
    =\frac{A_{k-1}}{A_k}x_{k-1}
     +\frac{\lambda_k}{A_k}y_k
    =\frac{1}{A_k}\sum_{i=1}^k\lambda_i y_i,
\]
and evaluate the new gradient at the averaged iterate $x_k$,
rather than at $y_k$, setting $g_k=\nabla f(x_k)$.

\subsection{Aggregated affine lower model for adaptive regularization}\label{subsec:agg-lower-model}

We now motivate our adaptive strategy for updating
the regularization parameters based on an accumulated affine lower model. Following the
estimate-sequence and approximate-duality-gap techniques~\cite{wibisono2015accelerated,
diakonikolas2019approximate}, we aggregate supporting affine functions and add a regularization term. 

For any $x\in\R^n$, we define the aggregated affine lower model $\varphi_k$ by
\begin{align}
    &\varphi_k(x)
    :=
    c_k+\langle s_k,x-x_0\rangle
    =
    \sum_{i=1}^k\lambda_i
    (f_i+\langle g_i,x-x_i\rangle)
    \le A_kf(x)\label{ineq:lower-ppty}\\
    &\text{with}\quad c_k:=
    \sum_{i=1}^k\lambda_i(f_i+\langle g_i,x_0-x_i\rangle),\nonumber
\end{align}
where the inequality is due to the convexity of $f$. This model collects the first-order information obtained
along the optimization trajectory. It requires only the vector $s_k$ and the scalar $c_k$, both of which can be updated recursively.

Directly minimizing $\varphi_k$ over $\mathbb{R}^n$ does not
generally yield a useful comparison value. We therefore
define the regularized estimating function for any $\gamma>0$ and its minimum
value by
\begin{align}\label{def:estimate-function-value}
    \psi_k(x;\gamma)
    :=\varphi_k(x)+\frac{\gamma}{2}\|x-x_0\|_p^2,\qquad
    m_k(\gamma)
    :=\min_{x\in\mathbb{R}^n}\psi_k(x;\gamma)
    =c_k-\frac{\|s_k\|_q^2}{2\gamma},
\end{align}
where the last equality follows from Fenchel conjugacy.
For any minimizer $x_*$ of $f$, \eqref{ineq:lower-ppty}
then implies
\begin{align}\label{mk-upbd-psi}
    m_k(\gamma)
    \le
    \psi_k(x_*;\gamma)
    \le
    A_kf_*+\frac{\gamma}{2}\|x_*-x_0\|_p^2.
\end{align}
Although $m_k(\gamma)$ is not necessarily a lower bound on
$A_kf_*$, this inequality provides the comparison needed
to control the objective gap.

Motivated by \eqref{mk-upbd-psi}, we seek a regularization
parameter satisfying
\begin{align}\label{Akfk-upbd-gamma-lineaer}
    A_kf_k
    \le
    m_k(\gamma)+\rho_k^2\gamma,
\end{align}
where $\rho_k>0$ is a scale to be specified later.
The additional term $\rho_k^2\gamma$ is chosen to scale
linearly with $\gamma$, matching the regularization term
in \eqref{mk-upbd-psi}. %, so that the resulting objective-gap bound is proportional to $\gamma$. 
The condition \eqref{Akfk-upbd-gamma-lineaer} is computable without knowing $f_*$.
Combining it with \eqref{mk-upbd-psi} yields
\begin{align}\label{ineq:upbd-fk-f*-gamma}
    f_k-f_*
    \le
    \frac{\gamma}{A_k}
    \Big(
        \rho_k^2+\frac12\|x_*-x_0\|_p^2
    \Big).
\end{align}
For a fixed $\rho_k$, the resulting bound is proportional
to $\gamma$. %, motivating the choice of the smallest admissible regularization parameter.

Given $\gamma_k>0$ and a safeguard parameter $\theta\in(0,1)$,
we select the smallest $\gamma\ge\theta\gamma_k$ satisfying
\eqref{Akfk-upbd-gamma-lineaer}:
\begin{align}\label{gammak+1-upper-lower}
    \gamma_{k+1}
    :=
    \min\left\{
        \gamma\ge\theta\gamma_k:
        A_kf_k\le m_k(\gamma)+\rho_k^2\gamma
    \right\}.
\end{align}
The safeguard limits the decrease in the regularization
parameter, while the estimating-function comparison
\eqref{Akfk-upbd-gamma-lineaer} determines whether a value
larger than $\theta\gamma_k$ is needed.
Denote $\chi_k:=A_kf_k-c_k$.
Using the expression for $m_k(\gamma)$ in \eqref{def:estimate-function-value}, we can rewrite
the update \eqref{gammak+1-upper-lower} explicitly as
\[
    \gamma_{k+1}
    =
    \max\left\{
        \theta\gamma_k,\,
        \frac{
            \chi_k+\sqrt{\chi_k^2+2\rho_k^2\|s_k\|_q^2}
        }{2\rho_k^2}
    \right\}.
\]
Since \eqref{Akfk-upbd-gamma-lineaer} holds with
$\gamma=\gamma_{k+1}$, it follows from
\eqref{ineq:upbd-fk-f*-gamma} that
\[
    f_k-f_*
    \le
    \frac{\gamma_{k+1}}{A_k}
    \Bigl(
        \rho_k^2+\frac12\|x_*-x_0\|_p^2
    \Bigr).
\]

\subsection{Bounded Lyapunov-potential certificate}\label{subsec:blpc}

We now use the estimating-function comparison in
Section~\ref{subsec:agg-lower-model} to construct a natural
potential sequence $\{\mathcal{P}_k\}$ for the convergence
analysis. Specifically, for all $k\ge1$, define
\[
    \mathcal{P}_k
    :=
    A_kf_k-m_k(\gamma_{k+1})
    =
    \chi_k+\frac{\|s_k\|_q^2}{2\gamma_{k+1}},
\]
where the last equality follows from the definition of
$\chi_k$ and \eqref{def:estimate-function-value}.
The connection between this potential and the objective
gap follows from \eqref{mk-upbd-psi}, which gives
\[
    \mathcal{P}_k
    \ge
    A_k(f_k-f_*)
    -\frac{\gamma_{k+1}}{2}\|x_*-x_0\|_p^2.
\]
Thus, controlling $\mathcal{P}_k$ relative to
$\gamma_{k+1}$ provides a bound on the objective gap. In addition, the condition \eqref{Akfk-upbd-gamma-lineaer}, with
$\gamma=\gamma_{k+1}$, ensures that
\[
    \mathcal{P}_k
    =
    \chi_k+\frac{\|s_k\|_q^2}{2\gamma_{k+1}}
    \le \rho_k^2\gamma_{k+1}.
\]
This bounded potential condition also provides a basis for establishing
the boundedness of the iterates $\{x_k\}$ and $\{y_k\}$. A detailed statement
will be given in Theorem~\ref{thm:bound-iterates}.

\subsection{The $\text{A}^2$ODA method and its convergence rate}\label{subsec:complete-alg}

We now summarize the $\text{A}^2$ODA method in
Algorithm~\ref{alg:pODA}, which combines the updates described above.

\begin{breakablealgorithm}%[!htpb]
\caption{$\text{A}^2$ODA: An \textbf{A}daptive \textbf{A}ccelerated \textbf{O}ptimistic \textbf{D}ual \textbf{A}veraging method}
\label{alg:pODA}
\begin{algorithmic}
\State\textbf{Input:} Initial iterate $x_0\in\R^n$ satisfying $\nabla f(x_0)\neq0$, a conjugate pair $(p,q)$ with $q\in[1,2]$ and $1/p+1/q = 1$, averaging scale $m\ge 2$, initial regularization parameter $\gamma_1>0$, and scaling parameters $\alpha,\theta$ satisfying $\theta\in(0,1)$ and $\alpha\!>\!2(1+e+e/\theta)$.
\State Set $f_0=f(x_0)$, $g_0=\nabla f(x_0)$, $s_0=0$, $c_0=0$, $A_0=0$, and {$\rho_0=0$}.
\For{$k=1,2,\ldots$}
    \State Compute the averaging weights:
    \begin{align}\label{def:Ak-lambdak}
        \lambda_k = \left( 1+ \frac{k-1}{m-1} \right)^{m-1}, \qquad A_k = A_{k-1} + \lambda_k.
    \end{align}
    \State Perform an optimistic dual averaging step:
    \begin{align}\label{alg1:AODA-steps}
        \Tilde{s}_k = s_{k-1} + \lambda_k\, g_{k-1},\qquad
        y_k = x_0 + \frac{1}{\gamma_k}\,J_q(-\Tilde{s}_k),\qquad
        x_k = \frac{A_{k-1}}{A_k}\,x_{k-1} + \frac{\lambda_k}{A_k}\, y_k.
    \end{align}
    \State Evaluate $f_k=f(x_k)$, $g_k=\nabla f(x_k)$, and update the safeguard radius:
    \begin{equation} \label{eq:update-rho}
    \rho_k=\max\bigg\{\rho_{k-1},\frac{\|y_k-x_0\|_p}{\alpha}\bigg\}.    
    \end{equation}
    \State Update the accumulated gradient and affine model:
    \begin{align}\label{def:sk-ck-chik}
        s_k = s_{k-1} + \lambda_k\, g_k,\qquad c_k = c_{k-1} + \lambda_k\,(f_k + \langle g_k, x_0-x_k\rangle),\qquad\chi_k = A_k\, f_k - c_k.
    \end{align}
    \State Update the regularization parameter:
    \begin{align}
        \gamma_{k+1} =\max\left\{\theta\gamma_k, \frac{\chi_k + \sqrt{\chi_k^2 + 2({\rho_k}\|s_k\|_q)^2}}{2{\rho_k^2}} \right\}.\label{def:betak+1}%\\
        %\gamma_{k+1} & = \frac{\gamma_{k+1}^{p-1}}{\max\{\|s_k\|_q, \|s_k + (k+1) g_k\|_q\}^{p-2}}.\label{def:gammak+1}
    \end{align}
\EndFor
\end{algorithmic}
\end{breakablealgorithm}

We now make a few remarks about Algorithm \ref{alg:pODA}.
\begin{itemize}[leftmargin=*]
\item The parameter $m\ge2$ controls the growth of the averaging
weights: for fixed $m$, $\lambda_k=\Theta(k^{m-1})$ and
$A_k=\Theta(k^m)$. In particular, $m=2$ gives $\lambda_k=k$
and $A_k=k(k+1)/2$, recovering the linearly increasing
weights used in accelerated optimistic dual
averaging~\cite{joulani2020averaging}.
Allowing $m>2$ provides flexibility to exploit the diminishing
gap-dependent curvature under generalized local smoothness,
potentially yielding faster decay of the corresponding
terms in the convergence bound. 

The benefit of increasing $m$ depends on its
interaction with the adaptive regularization. Indeed,
when the sequence of radius scales $\{\rho_k\}$ is bounded,
the objective-gap bound derived above is governed by
$\gamma_{k+1}/A_k$. Thus, faster growth of $A_k$ can translate into faster convergence only if it is not offset by the growth of $\gamma_{k+1}$. The convergence analysis
below quantifies this interaction.

\item The scale $\rho_k$ controls the tolerance in the estimating-function comparison (see \eqref{Akfk-upbd-gamma-lineaer}) and therefore affects the update of the regularization parameter $\gamma_{k+1}$. We choose
\[
    \rho_k=\frac{1}{\alpha}
    \max_{1\le i\le k}\{\|y_i-x_0\|_p\}
\]
so that this scale adapts to the observed size of the optimization trajectory, without requiring prior knowledge of the distance from $x_0$ to a minimizer. 

This choice balances two effects: a larger $\rho_k$ relaxes the estimating-function comparison and may allow a smaller regularization parameter, but it also increases the radius-dependent term in the objective-gap bound. This adaptive update is designed to achieve a favorable dependence on problem constants while using only quantities already available from the iterates. A simpler alternative is to use a fixed scale $\rho_k\equiv\rho>0$. 

\item The initial regularization parameter $\gamma_1$ can be selected through an optional backtracking line search, similar to that in~\cite{li2026simple}. We present its procedure as follows:

\smallskip
\refstepcounter{calibration}\label{alg:gamma1-calibration}
\begin{tcolorbox}[
    enhanced,
    breakable,
    colback=teal!2,
    colframe=teal!45!black,
    boxrule=1pt,
    arc=1mm,
    left=3mm,
    right=3mm,
    top=2.5mm,
    bottom=2.5mm,
    before skip=8pt,
    after skip=8pt
]
\noindent
{\textbf{Calibration~\thecalibration} (Selection of $\gamma_1$):} 
Given $g_0\ne0$ and an initial $\gamma>0$, define
$x_1(\gamma)=x_0+\gamma^{-1}J_q(-g_0)$.
A trial $\gamma$ is accepted if
\[
    f(x_1(\gamma))
    \le f(x_0)-\frac{\|g_0\|_q^2}{2\gamma}.
\]
Test the initial $\gamma$ and proceed as follows:
\begin{enumerate}[
    leftmargin=*,
    nosep,
    topsep=3pt,
    itemsep=4pt
]
\item \underline{\textit{Initial trial rejected.}}
   Repeatedly set $\gamma\leftarrow2\gamma$ and test the new $\gamma$ until it is accepted.
\item \underline{\textit{Initial trial accepted.}}
    Repeatedly test $\gamma/2$. If accepted,
    set $\gamma\leftarrow\gamma/2$;
    otherwise, stop.
\end{enumerate}
\smallskip
\noindent
\textbf{Return} $\gamma_1=\gamma$ and $x_1=x_1(\gamma_1)$.
\end{tcolorbox}
\smallskip

This procedure terminates after finitely many trials: differentiability of $f$ at $x_0$ ensures acceptance for all sufficiently large $\gamma$, while $f\ge f_*$ and $g_0\ne0$ ensure rejection for all sufficiently small $\gamma$.
Let $L_p$ and $H_p$ be the corresponding curvature moduli on a compact convex set $\mathcal{S}_p$ containing the iterates and a minimizer $x_*$, where $\mathcal{S}_p$ will be specified in Theorem \ref{thm:bound-iterates} below.
Under Assumption~\ref{asp:cvx}, the calibrated parameter $\gamma_1$ satisfies
\begin{equation} \label{pro:gamma1}
\frac{\|g_0\|_q^2}{2\Delta(x_0)}\le\gamma_1<2\bigl(L_p+H_p\Delta(x_0)^\nu\bigr).
\end{equation}
A proof of this bound is provided in Appendix~\ref{app:calibration-bound}.
The calibration requires no prior knowledge of the smoothness constants and improves the dependence of the convergence rate on $L_p$ and $H_p$, as detailed in Remark~\ref{re:calibrated-convergence-rate}. Since the line search is performed only during initialization, subsequent iterations remain line-search-free.
\end{itemize}

We conclude this section by presenting the convergence rate of the $\text{A}^2$ODA method under the generalized local smoothness condition in Assumption~\ref{asp:cvx}(b). A formal statement is provided in Theorem~\ref{thm:convex-convergence-rate} below. 
For any fixed $m\ge2$ and $\nu\in[0,1)$, with a calibration of $\gamma_1$, we establish
\begin{equation} \label{re:clear-rate}
\begin{aligned}
f(x_T)-f_*
=\mathcal{O}\Bigg(\frac{L_p\|x_0-x_*\|_p^2
+(H_p\|x_0-x_*\|_p^2)^{\frac{1}{1-\nu}}}{T^m} +\frac{L_p\|x_0-x_*\|_p^2}{T^{1+2/p}}
+\frac{(H_p\|x_0-x_*\|_p^2)^{\frac{1}{1-\nu}}}
{T^{\min\left\{m,\frac{1+2/p}{1-\nu}\right\}}}
\Bigg),
\end{aligned}
\end{equation}
with an additional logarithmic factor $(\log T)^{(1-2/p)/(1-\nu)}$ in the last term when $m=(1+2/p)/(1-\nu)$. The three terms reflect the effects of initialization, constant curvature, and gap-dependent curvature, respectively. In particular, when $\nu>0$, the gap-dependent curvature term can decay faster than the constant-curvature term. For $L_p>0$, the overall rate remains $\mathcal{O}(T^{-(1+2/p)})$, which matches the lower bound on the convergence-rate exponent in \cite{diakonikolas2020lower}, up to dimension-dependent factors in the high-dimensional oracle regime. When $L_p=0$ and $\nu>0$, a suitable choice of $m$ yields
a faster convergence rate. These results highlight the advantages of exploiting generalized local smoothness through adaptive regularization and weighted averaging.

\section{Convergence analysis}\label{sec:proof}

In this section, we establish convergence rate guarantees for the $\text{A}^2$ODA method by analyzing the potential sequence
\begin{align}\label{def:Lyapunov}
    \mathcal{P}_0=0,\qquad \mathcal{P}_k
    :=
    \chi_k+\frac{1}{\gamma_{k+1}}V_q(-s_k)
    \qquad \forall k\geq 1,
\end{align}
where $\{\chi_k\}$, $\{s_k\}$, and $\{\gamma_k\}$ are generated by the algorithm.

The following theorem establishes upper bounds on $\{\mathcal{P}_k\}$ and the boundedness of the iterates generated by Algorithm~\ref{alg:pODA}.

\begin{theorem}\label{thm:bound-iterates}
Suppose that Assumption~\ref{asp:cvx}(a) holds. Let
$\{(x_k,y_k,\gamma_k,\rho_k)\}$ be generated by Algorithm \ref{alg:pODA} 
% with $q\in[1,2]$, $\theta\in(0,1)$, \textcolor{blue}{$\gamma_1>0$, and $c>2(1+e+e/\theta)$}, let $p=q/(q-1)$, 
and let
$\{\mathcal{P}_k\}$ be defined in \eqref{def:Lyapunov}. {Fix any minimizer $x_*$ of $f$, and denote 
\begin{equation} \label{eq:bar-rho}
C_\theta=2(1+e+e/\theta), \quad \text{and}\quad \bar\rho=\max\Big\{\|g_0\|_q/(\alpha\gamma_1),\,C_\theta\|x_0-x_*\|_p/\sqrt{\alpha^2-C_\theta^2}\Big\}.
\end{equation}
Then, $0<\rho_k\le\bar\rho$ for $k\ge1$}, and for any $k\ge0$, 
\begin{align}
    \mathcal{P}_k\le {\rho_k^2}\gamma_{k+1}.
    \label{eq:potential-budget-bound}
\end{align}
Moreover, it holds that {$\|x_k-x_0\|_p,\|y_k-x_0\|_p\le \alpha\rho_k$} for any $k\ge1$, and $\|x_*-x_0\|_p\le\alpha\bar\rho$. Consequently, for all $k\ge1$,
\begin{align}
    x_k,y_k,x_*\in\mathcal S_{p}:=
    \left\{
        x\in\R^n:
        \|x-x_0\|_p
        \le \alpha\bar\rho
    \right\}.
    \label{def:bounded-set}
\end{align}
\end{theorem}

\begin{proof}
We first prove \eqref{eq:potential-budget-bound}. For $k=0$,
it follows from $\mathcal{P}_0=0$ and $\gamma_1>0$.
Using the update rule of $\rho_k$ given in \eqref{eq:update-rho}, we have $\rho_1=\|g_0\|_q/(\alpha\gamma_1)>0$ and $\{\rho_k\}_{k\geq 1}$ is a nondecreasing sequence, which implies that $\rho_k>0$ for all $k\ge1$.
Fix any $k\ge1$, by the definitions of $V_q$ and $\mathcal{P}_k$ in \eqref{def:V} and \eqref{def:Lyapunov}, one has
\begin{align}
    \mathcal{P}_k
    =
    \chi_k+\frac{1}{\gamma_{k+1}}V_q(-s_k)
    =
    \chi_k+\frac{\|s_k\|_q^2}{2\gamma_{k+1}}.
    \label{eq:quadratic-potential-expression}
\end{align}
Observe from the update of $\gamma_{k+1}$ in \eqref{def:betak+1} that
\[
    \gamma_{k+1}
    \ge
    \frac{
        \chi_k+
        \sqrt{\chi_k^2+2{\rho_k^2}\|s_k\|_q^2}
    }{2{\rho_k^2}}\qquad\Longrightarrow\qquad {\rho_k^2}\gamma_{k+1}^2
    -\chi_k\gamma_{k+1}
    -\frac12\|s_k\|_q^2
    \ge0.
\]
Dividing this inequality by $\gamma_{k+1}>0$ and using
\eqref{eq:quadratic-potential-expression}, we obtain
\[
    \mathcal{P}_k
    =
    \chi_k+\frac{\|s_k\|_q^2}{2\gamma_{k+1}}
    \le
    {\rho_k^2}\gamma_{k+1}.
\]
This proves \eqref{eq:potential-budget-bound}.
We next prove the boundedness of the iterates. Recall from \eqref{def:sk-ck-chik} and \eqref{def:Ak-lambdak} that
\begin{align}\label{eq:sum-express-sk-ck-Ak}
s_k=\sum_{i=1}^k \lambda_i\,g_i,\qquad
c_k=\sum_{i=1}^k \lambda_i\,(f_i+\langle g_i,x_0-x_i\rangle),\qquad A_k = \sum_{i=1}^k \lambda_i,     
\end{align}
hold for any $k\geq 1$. Using these identities, we obtain
\begin{align*}
    \mathcal{P}_k
    &\overset{\eqref{eq:quadratic-potential-expression}}{=}
    \chi_k+\frac{\|s_k\|_q^2}{2\gamma_{k+1}}
    \overset{\eqref{def:sk-ck-chik}, \eqref{def:Ak-lambdak}}{=}
    A_kf_k-c_k+\frac{\|s_k\|_q^2}{2\gamma_{k+1}} =
    A_k(f_k-f_*)
    +\sum_{i=1}^k \lambda_i\,D_f(x_*,x_i)
    +\langle s_k,x_*-x_0\rangle
    +\frac{\|s_k\|_q^2}{2\gamma_{k+1}} \\
    &\ge
    \frac{\|s_k\|_q^2}{2\gamma_{k+1}}
    +\langle s_k,x_*-x_0\rangle
    \ge
    \frac{\|s_k\|_q^2}{2\gamma_{k+1}}
    -\|s_k\|_q\|x_*-x_0\|_p \ge
    \frac{\|s_k\|_q^2}{4\gamma_{k+1}}
    -\gamma_{k+1}\|x_*-x_0\|_p^2,
\end{align*}
where the third equality follows from \eqref{eq:sum-express-sk-ck-Ak} and the definition of $D_f$, the first inequality is due to $f_k\ge f_*$ and the convexity of $f$, and the last two
inequalities follow from H\"older's and Young's inequalities,
respectively. Combining this inequality with $\mathcal{P}_k\le {\rho_k^2}\gamma_{k+1}$ yields
\begin{equation}
\label{eq:scaled-gradient-boundedness}
    \frac{\|s_k\|_q}{\gamma_{k+1}}
    \le
    2\sqrt{\|x_*-x_0\|_p^2+{\rho_k^2}}.
\end{equation}
This also holds for $k=0$ because $s_0=0$. By \eqref{ppty:J} and the update of $y_k$ in \eqref{alg1:AODA-steps}, one has for any $k\ge1$ that
\begin{align}\label{eq:diff-yk-x0}
    \|y_k-x_0\|_p
    =
    \frac{1}{\gamma_k}\|J_q(-\tilde s_k)\|_p
    =
    \frac{\|\tilde s_k\|_q}{\gamma_k}.
\end{align}
For $k=1$, the update rule gives
$\|y_1-x_0\|_p=\|g_0\|_q/\gamma_1=\alpha\rho_1$.
For $k\ge2$, combining $s_{k-1}-s_{k-2}=\lambda_{k-1}g_{k-1}$ and
$\tilde s_k-s_{k-1}=\lambda_kg_{k-1}$ gives
\begin{align}\label{def:t-sk-idt}
    \tilde s_k
    =
    \Big(1+\frac{\lambda_k}{\lambda_{k-1}}\Big)s_{k-1}
    -\frac{\lambda_k}{\lambda_{k-1}}s_{k-2}.
\end{align}
By Lemma~\ref{lem:weights-estimates}, we have $\lambda_k/\lambda_{k-1}\le e$.
Combining with $\gamma_k\ge\theta\gamma_{k-1}$, and $\rho_{k-2}\le\rho_{k-1}$ gives
\begin{align*}
\|y_k-x_0\|_p
& \overset{\eqref{eq:diff-yk-x0}}{=} \frac{\|\tilde{s}_k\|_q}{\gamma_k}\le
\Big(1+\frac{\lambda_k}{\lambda_{k-1}}\Big)
\frac{\|s_{k-1}\|_q}{\gamma_k}
+\frac{\lambda_k}{\lambda_{k-1}}
\frac{\|s_{k-2}\|_q}{\gamma_k}\leq
(1+e)\frac{\|s_{k-1}\|_q}{\gamma_k}
+\frac{e}{\theta}
\frac{\|s_{k-2}\|_q}{\gamma_{k-1}}\\
&\overset{\eqref{eq:scaled-gradient-boundedness}}{\le}
2(1+e)\sqrt{\|x_0-x_*\|_p^2+\rho_{k-1}^2}
+\frac{2e}{\theta}
\sqrt{\|x_0-x_*\|_p^2+\rho_{k-2}^2}\\
&\le
2\Big(1+e+\frac{e}{\theta}\Big)
\sqrt{\|x_0-x_*\|_p^2+\rho_{k-1}^2}\overset{\eqref{eq:bar-rho}}{=} C_\theta
\sqrt{\|x_0-x_*\|_p^2+\rho_{k-1}^2}.
\end{align*}
Consequently, by \eqref{eq:update-rho} and the preceding
bound on $\|y_k-x_0\|_p$, for every $k\ge2$, we have
\[
\rho_k\le
\max\left\{
    \rho_{k-1},
    \frac{C_\theta}{\alpha}
    \sqrt{\|x_0-x_*\|_p^2+\rho_{k-1}^2}
\right\}.
\]
Since in Algorithm \ref{alg:pODA}, we set $\alpha>C_\theta$, the definition of $\bar\rho$ in \eqref{eq:bar-rho} gives
$\rho_1=\|g_0\|_q/(\alpha\gamma_1)\le\bar\rho$ and
\[
C_\theta^2\|x_0-x_*\|_p^2
\le(\alpha^2-C_\theta^2)\bar\rho^2   \qquad\Longrightarrow\qquad \frac{C_\theta}{\alpha}
\sqrt{\|x_0-x_*\|_p^2+\bar\rho^2}
\le\bar\rho.
\]
We now prove $\rho_k\le\bar\rho$ by induction.
The case $k=1$ has already been established.
Suppose that $\rho_{k-1}\le\bar\rho$ for some $k\ge2$, we have
\[
\begin{aligned}
\rho_k
\le
\max\left\{
    \bar\rho,
    \frac{C_\theta}{\alpha}
    \sqrt{\|x_0-x_*\|_p^2+\bar\rho^2}
\right\}=\bar\rho.
\end{aligned}
\]
Thus, we have $\rho_k\le\bar\rho$ for all $k\ge1$.
It remains to bound the iterates.
By \eqref{eq:update-rho},
$\|y_k-x_0\|_p\le \alpha \rho_k$ for every $k\ge1$.
To obtain the corresponding bound for $x_k$,
note first that $\|x_0-x_0\|_p=0=\alpha\rho_0$.
Suppose that
$\|x_{k-1}-x_0\|_p\le \alpha\rho_{k-1}$.
Using the update of $x_k$ in
\eqref{alg1:AODA-steps},
$A_k=A_{k-1}+\lambda_k$, and
$\rho_{k-1}\le\rho_k$, we obtain
\begin{align*}
\|x_k-x_0\|_p
&=
\left\|
    \frac{A_{k-1}}{A_k}(x_{k-1}-x_0)
    +\frac{\lambda_k}{A_k}(y_k-x_0)
\right\|_p\le
\frac{A_{k-1}}{A_k}
\|x_{k-1}-x_0\|_p
+\frac{\lambda_k}{A_k}
\|y_k-x_0\|_p \\
&\le
\alpha\left(
    \frac{A_{k-1}}{A_k}\rho_{k-1}
    +\frac{\lambda_k}{A_k}\rho_k
\right)
\le
\alpha\frac{A_{k-1}+\lambda_k}{A_k}\rho_k=\alpha\rho_k.
\end{align*}
Combining these bounds with $\rho_k\le\bar\rho$ yields
$x_k,y_k\in\mathcal S_p$.
Moreover, since $\alpha>C_\theta>1$, the definition of $\bar\rho$ in \eqref{eq:bar-rho} gives
\[
\alpha\bar\rho
\ge
\frac{\alpha C_\theta}
{\sqrt{\alpha^2-C_\theta^2}}
\|x_0-x_*\|_p
>
\|x_0-x_*\|_p,
\]
which also implies $x_*\in\mathcal S_p$.
This completes the proof.
\end{proof}

\begin{remark}[Upper bounds on $\bar\rho$]
\label{re:barrho}
We now derive upper bounds on $\bar\rho$ for the cases with and without calibration of $\gamma_1$. The calibration procedure is described in Calibration \ref{alg:gamma1-calibration}.
% .~\ref{subsec:complete-alg}.
%We first relate $\bar\rho$ to the initial distance for arbitrary $\gamma_1>0$ and then discuss an optional startup line search for selecting $\gamma_1$.
\begin{itemize}[leftmargin=1.2em]
\item[(i)] \underline{Without calibration of $\gamma_1$}.
Since $\nabla f(x_*)=0$, \eqref{ineq:general-smooth} gives
$\|g_0\|_q\le
(L_p+H_p\Delta(x_0)^\nu)\|x_0-x_*\|_p$.
Hence, the definition of $\bar\rho$ yields
\[
\bar\rho\le
\max\left\{
\frac{L_p+H_p\Delta(x_0)^\nu}{\alpha\gamma_1},
\frac{C_\theta}{\sqrt{\alpha^2-C_\theta^2}}
\right\}\|x_0-x_*\|_p.
\]

\item[(ii)]\underline{With calibration of $\gamma_1$}.
Since $f(x_1)\ge f_*$, the
sufficient-decrease condition and convexity imply
\[
    \frac{\|g_0\|_q^2}{2\gamma_1}
    \le \Delta(x_0)
    \le \langle g_0,x_0-x_*\rangle
    \le \|g_0\|_q\|x_0-x_*\|_p.
\]
Consequently, $\|g_0\|_q/(\alpha\gamma_1)\le(2/\alpha)\|x_0-x_*\|_p$. Hence, the definition of $\bar\rho$ yields
\[
    \bar\rho
    \le
    \max\left\{
        \frac{2}{\alpha},
        \frac{C_\theta}{\sqrt{\alpha^2-C_\theta^2}}
    \right\}
    \|x_0-x_*\|_p.
\]
\end{itemize}
Therefore, calibrating $\gamma_1$ eliminates the explicit dependence on $L_p$ and $H_p$ in the upper bound for $\bar\rho$, leading to improved dependence on these constants in the convergence rate. %This conclusion holds for all $\nu\in[0,1)$.
\end{remark}

Using Assumption~\ref{asp:cvx}(b) and the compactness and convexity
of $\mathcal S_p$, there exist constants $L_p(\mathcal S_p), H_p(\mathcal S_p)\ge0$, denoted simply by $L_p, H_p$, such that
\begin{align}\label{ineq:bar-Lip-p}
    \|\nabla f(y)-\nabla f(x)\|_q^2
    \le 2(L_p + H_p\max\{\Delta(x),\Delta(y)\}^\nu)D_f(x,y)
    \qquad\forall x,y\in\mathcal S_p.
\end{align}
For each $k\ge2$, define the auxiliary curvature estimate by
\begin{align}\label{def:ellk}
    \ell_k
    :=
    \begin{cases}
        h_k\dfrac{\|g_k-g_{k-1}\|_q^2}
                  {2D_f(x_{k-1},x_k)},
        & D_f(x_{k-1},x_k)>0,\\[1.5mm]
        0,
        & D_f(x_{k-1},x_k)=0,
    \end{cases}
    \qquad
    h_k:=\frac{\lambda_k^2}{2A_{k-1}}.
\end{align}
For each $k\ge1$, define the local curvature bound
\begin{align}\label{def:Mk}
    M_k:=L_p+H_p\max\{\Delta(x_{k-1}),\Delta(x_k)\}^{\nu}.
\end{align}
Since $x_{k-1},x_k\in\mathcal S_p$, applying~\eqref{ineq:bar-Lip-p} yields $\|g_k-g_{k-1}\|_q^2\le 2M_kD_f(x_{k-1},x_k)$. Consequently,
\begin{align}\label{ineq:upbd-local-estimate}
    \ell_k
    \le h_k M_k
    \qquad\forall k\ge2.
\end{align}

Our convergence analysis relies on comparisons with reference
potentials. Specifically, we define the reference regularization
scales $\{\eta_{k}\}$ and the reference potentials
$\{\mathcal Q_k\}$ by
\begin{align}\label{def:eta-Q-analysis}
    \eta_{k+1}
    :=\|[\gamma_k,r_a\ell_k]\|_{r_q},
    \qquad
    \mathcal Q_k
    :=
    \chi_k+\frac{1}{\eta_{k+1}}V_q(-s_k)
    \qquad \forall k\ge2,
\end{align}
where $r_a>8$ is fixed and $r_q$ is defined in \eqref{def:rq}. For notational convenience, we also define the auxiliary sequences denoting displacements from $x_0$:
\begin{align}\label{def:zk-zk-}
    z_k
    :=\frac{1}{\gamma_k}J_q(-\tilde s_k),
    \qquad
    z_k^-
    :=\frac{1}{\gamma_k}J_q(-s_{k-1})
    \qquad \forall k\ge1.
\end{align}

\begin{lemma}
    Suppose that Assumption \ref{asp:cvx} holds. Let $\{(x_k,s_k,\tilde{s}_k, \rho_k,\gamma_k)\}$ be generated by Algorithm \ref{alg:pODA}, let $\{\mathcal{P}_k\}$ be defined in \eqref{def:Lyapunov}, let $\{\eta_k\}$ and $\{\mathcal{Q}_k\}$ be defined in \eqref{def:eta-Q-analysis}, and let $\{(z_k,z_k^-)\}$ be defined in \eqref{def:zk-zk-}. Then, for any $k\ge2$,
    \begin{align}
    \mathcal{Q}_k - \mathcal{P}_{k-1} & =
    -\frac{1}{\gamma_k}V_q(-\tilde s_k)
    +\langle z_k,s_k-\tilde s_k\rangle + \frac{1}{\eta_{k+1}}V_{q}(-s_k)\nonumber \\
    & \qquad -A_{k-1}D_f(x_{k-1},x_k) -\gamma_kD_{d_p}^{-s_{k-1}/\gamma_k}(z_k,z_k^-).\label{eq:Lypnov-idtt}
    \end{align}
\end{lemma}

\begin{proof}
Fix any $k\ge2$. By \eqref{alg1:AODA-steps}, \eqref{def:Ak-lambdak}, and \eqref{def:zk-zk-}, one has that $A_{k-1}(x_k-x_{k-1})= \lambda_k(y_k-x_k)$ and $z_k=y_k-x_0$. Using these and the updates of $A_k$, $c_k$, and $\chi_k$ in \eqref{def:sk-ck-chik} and \eqref{def:Ak-lambdak}, we have
\begin{align*}
    \chi_k-\chi_{k-1}
    & =
    A_kf_k-A_{k-1}f_{k-1}
    - \lambda_k(f_k+\langle g_k,x_0-x_k\rangle)=
    A_{k-1}(f_k-f_{k-1})
    + \lambda_k\langle g_k,x_k-x_0\rangle\\
    & = A_{k-1}\langle g_k,x_k-x_{k-1}\rangle + \lambda_k\langle g_k,x_k-x_0\rangle - A_{k-1}D_f(x_{k-1},x_k)\\
    & = \lambda_k \langle g_k, y_k-x_0\rangle - A_{k-1}D_f(x_{k-1},x_k) = \lambda_k \langle g_k, z_k\rangle - A_{k-1}D_f(x_{k-1},x_k),
\end{align*}
where the third equality is due to the definition of $D_f$. By this and the definition of $\mathcal{P}_{k-1}$ in \eqref{def:Lyapunov}, one has
\begin{align}
    \chi_k-\mathcal{P}_{k-1}
    &=
    \lambda_k\langle\nabla f(x_k),z_k\rangle
    -A_{k-1}D_f(x_{k-1},x_k)-\frac{1}{\gamma_k}V_q(-s_{k-1}).
    \label{eq:Phi-increment-intermediate}
\end{align}
Since $V_q=d_p^*$, it follows from \eqref{def:J}, \eqref{def:zk-zk-}, and conjugacy that $-\tilde s_k/\gamma_k\in\partial d_p(z_k)$ and $-s_{k-1}/\gamma_k\in\partial d_p(z_k^-)$. By Fenchel-Young equality in \eqref{eq:FY}, one has
\begin{equation}\label{eq:FY-1}
    \gamma_k d_p(z_k)
    =-\langle\tilde s_k,z_k\rangle
      -\frac{1}{\gamma_k}V_q(-\tilde s_k),\qquad \gamma_k d_p(z_k^-)
    =-\langle s_{k-1},z_k^-\rangle
      -\frac{1}{\gamma_k}V_q(-s_{k-1}).
\end{equation}
Consequently,
\begin{align*}
    \gamma_k D_{d_p}^{-s_{k-1}/\gamma_k}(z_k,z_k^-)
    &=
    \gamma_k d_p(z_k)-\gamma_k d_p(z_k^-)
    +\langle s_{k-1},z_k-z_k^-\rangle\\
    &\overset{\eqref{eq:FY-1}}{=}
    -\langle\tilde s_k,z_k\rangle
    -\frac{1}{\gamma_k}V_q(-\tilde s_k)
    +\langle s_{k-1},z_k^-\rangle
    +\frac{1}{\gamma_k}V_q(-s_{k-1})
    +\langle s_{k-1},z_k-z_k^-\rangle\\
    &=
    \frac{1}{\gamma_k}
    \bigl(V_q(-s_{k-1})-V_q(-\tilde s_k)\bigr)
    -\langle\tilde s_k-s_{k-1},z_k\rangle \\
    & =
    \frac{1}{\gamma_k}
    \bigl(V_q(-s_{k-1})-V_q(-\tilde s_k)\bigr)
    -\lambda_k\langle g_{k-1},z_k\rangle,
\end{align*}
where the last equality follows from
$\tilde s_k-s_{k-1}=\lambda_k g_{k-1}$.
Substituting this into
\eqref{eq:Phi-increment-intermediate} and using
$s_k-\tilde s_k=\lambda_k(g_k-g_{k-1})$, we obtain
\begin{align*}
    \chi_k-\mathcal P_{k-1}
    &=
    -\frac{1}{\gamma_k}V_q(-\tilde s_k)
    +\lambda_k\langle g_k-g_{k-1},z_k\rangle
    -A_{k-1}D_f(x_{k-1},x_k)
    -\gamma_kD_{d_p}^{-s_{k-1}/\gamma_k}(z_k,z_k^-)\\
    &=
    -\frac{1}{\gamma_k}V_q(-\tilde s_k)
    +\langle z_k,s_k-\tilde s_k\rangle
    -A_{k-1}D_f(x_{k-1},x_k)
    -\gamma_kD_{d_p}^{-s_{k-1}/\gamma_k}(z_k,z_k^-).
\end{align*}
Combining this with the definition of $\mathcal Q_k$ in
\eqref{def:eta-Q-analysis} completes the proof.
\end{proof}

We define the threshold
\begin{align}\label{def:aq}
    a_q:=\frac12-\frac{2^{2/r_q}}{r_a}>0,
\end{align}
where $r_q$ is defined in \eqref{def:rq}. We partition the iteration indices $k\ge2$ as follows:
\begin{align}\label{def:set-A-B}
\mathcal A:=\{k\ge2:\ell_k\le a_q\gamma_k\},
\qquad
\mathcal B:=\{k\ge2:\ell_k>a_q\gamma_k\}.
\end{align}
This partition is based on a comparison between the regularization
parameters $\{\gamma_k\}$ used by Algorithm~\ref{alg:pODA} and the local
curvature estimates $\{\ell_k\}$.

The next lemma shows that, for $k\in\mathcal A$, the next
regularization parameter satisfies $\gamma_{k+1}\le\eta_{k+1}$,
which indicates that Algorithm~\ref{alg:pODA} takes a larger step than
the reference step.

\begin{lemma}
\label{lem:descent-iterations}
Suppose that Assumption~\ref{asp:cvx} holds. Let $\{(\rho_k,\gamma_k)\}$ be generated by Algorithm \ref{alg:pODA},
and let $\{\mathcal P_k\}$, $\{\eta_k\}$, $\{\mathcal Q_k\}$, and $\mathcal A$ be defined in
\eqref{def:Lyapunov}, \eqref{def:eta-Q-analysis}, and \eqref{def:set-A-B}, respectively. Then, we have $\mathcal Q_k\le \mathcal P_{k-1}\le \rho_k^2\eta_{k+1}$ and $\gamma_{k+1}\le \eta_{k+1}$ for all $k\in\mathcal{A}$.
\end{lemma}

\begin{proof}
Fix any $k\in\mathcal{A}$. If $\ell_k=0$, it follows from \eqref{def:ellk} and the cocoercivity of $f$ that $g_k=g_{k-1}$. Then, one has $s_k=\tilde{s}_k$ and $\eta_{k+1}=\gamma_k$. Using these and \eqref{eq:Lypnov-idtt}, we obtain that
\begin{align*}
    \mathcal{Q}_k - \mathcal{P}_{k-1} & = \frac{1}{\eta_{k+1}}V_{q}(-s_k)
    -\frac{1}{\gamma_k}V_{q}(-\tilde s_k)
    +\langle z_k,s_k-\tilde s_k\rangle-A_{k-1}D_f(x_{k-1},x_k) -{\gamma_k}D_{d_p}^{-s_{k-1}/\gamma_k}(z_k,z_k^-)\\
    & = -A_{k-1}D_f(x_{k-1},x_k)-{\gamma_k}D_{d_p}^{-s_{k-1}/\gamma_k}(z_k,z_k^-) \le 0. 
\end{align*}
This together with $\eta_{k+1}=\gamma_k$ proves the conclusion.

Now suppose $\ell_k>0$. We next prove $\mathcal{Q}_k\le\mathcal{P}_{k-1}$ by considering two separate cases.

{(i)} $q=2$. By \eqref{def:zk-zk-} and \eqref{def:rq}, one has that $z_k=-\tilde s_k/\gamma_k$ and $r_q=\infty$. It then follows from \eqref{def:eta-Q-analysis} that $\eta_{k+1}=\max\{\gamma_k,r_a\ell_k\}\ge\gamma_k$. Then, the Euclidean identity \eqref{eq:E-idt} gives
\begin{align*}
    &\quad\frac{1}{\eta_{k+1}}V_{2}(-s_k)
    -\frac{1}{\gamma_k}V_2(-\tilde s_k)
    +\langle z_k,s_k-\tilde s_k\rangle\\
    &=
    \frac{\|s_k\|_2^2}{2\eta_{k+1}}
    -\frac{\|\tilde s_k\|_2^2}{2\gamma_k}
    -\frac{\langle\tilde s_k,s_k-\tilde s_k\rangle}{\gamma_k} =
    \frac{\|s_k-\tilde s_k\|_2^2}{2\gamma_k}
    -\frac12
    \Big(\frac1{\gamma_k}-\frac1{\eta_{k+1}}\Big)
    \|s_k\|_2^2 \le
    \frac{\|s_k-\tilde s_k\|_2^2}{2\gamma_k},
\end{align*}
where the last inequality is due to $\eta_{k+1}\ge\gamma_k$. The definition of $\ell_k$ in \eqref{def:ellk}, together with
$s_k-\tilde s_k=\lambda_k(g_k-g_{k-1})$ (see \eqref{alg1:AODA-steps} and \eqref{def:sk-ck-chik}), yields
\[
    A_{k-1}D_f(x_{k-1},x_k)
    =\frac{\|s_k-\tilde s_k\|_2^2}{4\ell_k}.
\]
Substituting the above two relations into the potential identity \eqref{eq:Lypnov-idtt} gives
\begin{align*}
    \mathcal Q_k-\mathcal P_{k-1}
    \le
    \Big(
        \frac1{2\gamma_k}-\frac1{4\ell_k}
    \Big)
    \|s_k-\tilde s_k\|_2^2
    -{\gamma_k}D_{d_2}(z_k,z_k^-) \le
    -{\gamma_k}D_{d_2}(z_k,z_k^-)
    \le0,
\end{align*}
where the second inequality follows from
$k\in\mathcal A$ and
$\ell_k\le a_2\gamma_k\le\gamma_k/2$,
and the last inequality follows from the convexity of
$d_2$.

{(ii)} $q\in[1,2)$. Since $\ell_k>0$, the definition of $\eta_{k+1}$ in
\eqref{def:eta-Q-analysis} gives $\eta_{k+1}>\gamma_k$ and $(\eta_{k+1}^{r_q}-\gamma_k^{r_q})^{1/r_q}=r_a\ell_k$. Applying \eqref{ineq:V-expan} with
$u=-\tilde s_k$, $v=\tilde s_k-s_k$, and
$c=\gamma_k/\eta_{k+1}\in(0,1)$, and using
$z_k=J_q(-\tilde s_k)/\gamma_k$, we obtain
\begin{align}\label{eq:residual-term}
    \frac{1}{\eta_{k+1}}V_q(-s_k) &
    -\frac{1}{\gamma_k}V_q(-\tilde s_k)
    +\langle z_k,s_k-\tilde s_k\rangle=
    \frac{1}{\gamma_k}
    \bigg(
        \frac{\gamma_k}{\eta_{k+1}}V_q(-s_k)
        -V_q(-\tilde s_k)
        +\langle J_q(-\tilde s_k),s_k-\tilde s_k\rangle
    \bigg) \nonumber\\
    &\le
    \frac{1}{\gamma_k}
    \bigg(
        1+
        \frac{2^{2/r_q}\gamma_k}
        {(\eta_{k+1}^{r_q}-\gamma_k^{r_q})^{1/r_q}}
    \bigg)
    V_q(s_k-\tilde s_k) =
    \bigg(
        \frac{1}{2\gamma_k}
        +\frac{2^{2/r_q}}{2r_a\ell_k}
    \bigg)
    \|s_k-\tilde s_k\|_q^2.
\end{align}
The definition of $\ell_k$ in \eqref{def:ellk}, together with
$s_k-\tilde s_k=\lambda_k(g_k-g_{k-1})$,
yields
\begin{align}\label{eq:ADfxk}
    A_{k-1}D_f(x_{k-1},x_k)
    =\frac{\|s_k-\tilde s_k\|_q^2}{4\ell_k}.
\end{align}
Substituting \eqref{eq:residual-term} and \eqref{eq:ADfxk}
into the potential identity \eqref{eq:Lypnov-idtt} gives
\begin{align*}
    \mathcal Q_k-\mathcal P_{k-1}
    &\le
    \Big(
        \frac{1}{2\gamma_k}
        +\frac{2^{2/r_q}}{2r_a\ell_k} - \frac{1}{4\ell_k}
    \Big)
    \|s_k-\tilde s_k\|_q^2
    -{\gamma_k} D_{d_p}^{-s_{k-1}/\gamma_k}(z_k,z_k^-)\\
    &=
    \frac12
    \Big(
        \frac{1}{\gamma_k}-\frac{a_q}{\ell_k}
    \Big)
    \|s_k-\tilde s_k\|_q^2
    -{\gamma_k} D_{d_p}^{-s_{k-1}/\gamma_k}(z_k,z_k^-) \le
    -{\gamma_k} D_{d_p}^{-s_{k-1}/\gamma_k}(z_k,z_k^-)
    \le0,
\end{align*}
where the equality uses \eqref{def:aq}, the second inequality is due to $k\in\mathcal A$ and $\ell_k\le a_q\gamma_k$,
and the last inequality follows from the convexity of $d_p$.

Further, by the definition of $\eta_{k+1}$ in \eqref{def:eta-Q-analysis}, one has $\eta_{k+1}\ge\gamma_k$. Then, using \eqref{eq:potential-budget-bound} and the definition of $\mathcal{Q}_k$ in \eqref{def:eta-Q-analysis}, we have
\begin{align*}
\chi_k + \frac{\|s_k\|_q^2}{2\eta_{k+1}}= \mathcal{Q}_k \le \mathcal{P}_{k-1} \le {\rho_{k-1}^2}\gamma_k \le {\rho_k^2}\eta_{k+1}\qquad\Longrightarrow\qquad \eta_{k+1}\ge\frac{
        \chi_k+
        \sqrt{\chi_k^2+2{\rho_k^2}\|s_k\|_q^2}
    }{2\rho_k^2}.
\end{align*}
This along with \eqref{def:betak+1} and $\eta_{k+1}\ge\gamma_k\ge\theta\gamma_k$ implies $\eta_{k+1}\ge\gamma_{k+1}$, which completes the proof.
\end{proof}

The next lemma bounds $\gamma_{k+1}$ by the scaled local curvature estimate, with a multiplicative constant independent of the iteration index.

\begin{lemma}\label{lem:gamma-bound-B}
Suppose that Assumption \ref{asp:cvx} holds. Let $\{(\rho_k,\gamma_k)\}$ be generated by Algorithm \ref{alg:pODA}, 
let $L_p$ and $H_p$ satisfy \eqref{ineq:bar-Lip-p}, and let $a_q$ be defined in \eqref{def:aq}, and recall that $\alpha>2(1+e+e/\theta)$, where $\theta\in (0,1)$ is the scaling parameter specified in Algorithm~\ref{alg:pODA}.
Define 
    \begin{align}\label{def:GammaB}
        \Gamma_{\mathcal{B}} := \Big(4 + \frac{1}{a_q}\Big) (1 + \alpha )^2.
    \end{align}
    Then, $\gamma_{k+1}\le \Gamma_{\mathcal{B}}h_k M_k$ for all $k\in\mathcal{B}$, where $\{h_k\}$ and $\{M_k\}$ are defined in \eqref{def:ellk} and \eqref{def:Mk}, respectively.
\end{lemma}

\begin{proof}
Fix any $k\in\mathcal{B}$. Recall from \eqref{ineq:upbd-local-estimate} that $\ell_k\le h_k M_k$. By this and the definition of $\mathcal{B}$, one has
\begin{align}\label{ineq:upbd-gammak-caseB}
    \gamma_k < \frac{\ell_k}{a_q} \le \frac{h_k M_k}{a_q}.
\end{align}
Recall that the update of $\gamma_{k+1}$ has the form
\[
    \gamma_{k+1}
    =\max\left\{
        \theta\gamma_k,\,
        \frac{\chi_k+\sqrt{\chi_k^2+2({\rho_k}\|s_k\|_q)^2}}
             {2{\rho_k}^2}
    \right\},
\]
it remains to derive upper bounds on $\chi_k$ and $\|s_k\|_q$
to bound $\gamma_{k+1}$. By the updates of $x_k$, $\tilde s_k$, and $s_k$ in \eqref{alg1:AODA-steps} and \eqref{def:sk-ck-chik}, the bounds $\|y_k-x_0\|_p,\|x_k-x_0\|_p,\|x_{k-1}-x_0\|_p\le \alpha\rho_k$ from Theorem~\ref{thm:bound-iterates}, and \eqref{ineq:general-smooth}, one has 
\begin{align}\label{ineq:sk-tsk-B}
    \|s_k-\tilde s_k\|_q = \lambda_k\|g_k - g_{k-1}\|_q \le \lambda_kM_k\|x_k - x_{k-1}\|_p = \frac{\lambda_k^2M_k}{A_k}\|y_k - x_{k-1}\|_p \le 4h_kM_k \alpha\rho_k.
\end{align}
In addition, using the definition of $z_k$ in \eqref{def:zk-zk-} and the updates of $y_k$, we have
\begin{align}\label{ineq:zk-tsk-B}
    \|z_k\|_p = \|y_k-x_0\|_p \le {\alpha\rho_k}, \qquad \|\tilde s_k\|_q \le \gamma_k\|z_k\|_p \le \gamma_k {\alpha\rho_k}.
\end{align}
Combining the identity \eqref{eq:Lypnov-idtt} with the definition of $\mathcal{Q}_k$ in \eqref{def:eta-Q-analysis}, we obtain
\begin{align*}
    \chi_k & = \mathcal{P}_{k-1} - \frac{1}{\gamma_k}V_{q}(-\tilde s_k)
    +\langle z_k,s_k-\tilde s_k\rangle -A_{k-1}D_f(x_{k-1},x_k) - {\gamma_k} D_{d_p}^{-s_{k-1}/\gamma_k}(z_k,z_k^-)\\
    &\le {\rho_k}^2\gamma_k + \langle z_k,s_k-\tilde s_k\rangle \le {\rho_k}^2\gamma_k + \|z_k\|_p \|s_k-\tilde s_k\|_q \overset{\eqref{ineq:upbd-gammak-caseB}\eqref{ineq:sk-tsk-B}\eqref{ineq:zk-tsk-B}}{\le} \frac{{\rho_k}^2h_kM_k}{a_q} + 4h_kM_k\alpha^2\rho_k^2,
\end{align*}
where the first inequality is due to \eqref{eq:potential-budget-bound}, {$\rho_{k-1}\le\rho_k$}, $V_{q}(-\tilde s_k)\ge0$, and the convexity of $f$ and $d_p$. Moreover, it follows from \eqref{ineq:upbd-gammak-caseB}, \eqref{ineq:sk-tsk-B}, and \eqref{ineq:zk-tsk-B} that
\begin{align*}
    \|s_k\|_q \le \|\tilde s_k\|_q + \|s_k-\tilde s_k\|_q \le \gamma_k\alpha\rho_k+ 4h_kM_k\alpha\rho_k \le \frac{h_kM_k\alpha\rho_k}{a_q} + 4h_kM_k\alpha\rho_k.
\end{align*} 
Combining the upper bounds for $\chi_k$ and $\|s_k\|_q$, we obtain
\begin{align*}
    \frac{\chi_k+\sqrt{\chi_k^2+2({\rho_k}\|s_k\|_q)^2}}
         {2{\rho_k}^2} &\le
    \frac{\chi_k+|\chi_k|}{2{\rho_k}^2}
    +\frac{\|s_k\|_q}{\sqrt{2}{\rho_k}} \le
   {\left(a_q^{-1}+4\alpha^2+\frac \alpha{\sqrt2}(a_q^{-1}+4)\right)h_kM_k} \\
    &\le
    \Gamma_{\mathcal{B}}h_kM_k,
\end{align*}
where the first inequality uses $\sqrt{a^2+b^2}\le |a|+|b|$ for all $a,b\in\R$, and the last inequality follows from the direct calculation. This inequality, together with \eqref{def:GammaB}, \eqref{ineq:upbd-gammak-caseB}, and \eqref{def:betak+1} implies $\gamma_{k+1}\le\Gamma_{\mathcal{B}}h_kM_k$. Hence, the conclusion of this lemma holds.
\end{proof}

Next, we introduce the quantity
\begin{equation} \label{def:gamma0}
\Gamma_0=\max\{\gamma_1,\gamma_2\} 
\end{equation}
to separate the initialization from the subsequent curvature estimates, while retaining $\gamma_1>0$ as an arbitrary input.

\begin{lemma}\label{lem:initial-scale-bound}
Suppose that Assumption~\ref{asp:cvx} holds. Let $\{\gamma_k\}$ be generated by Algorithm~\ref{alg:pODA}, and let $\Gamma_0$ and $M_1$ be defined in \eqref{def:gamma0} and \eqref{def:Mk}, respectively. Then, $\Gamma_0\le\alpha^2M_{1}
+{\alpha(\gamma_1+M_{1})}/{\sqrt2}$.
\end{lemma}
\begin{proof}
Since $\lambda_1=A_1=1$, the first iteration gives $x_1=y_1$, $s_1=g_1$, and $\chi_1=\langle g_1,x_1-x_0\rangle$. Moreover, it is easy to verify that $\|x_1-x_0\|_p=\alpha\rho_1$, $\|g_0\|_q=\gamma_1 \alpha\rho_1$, and $\langle g_0,x_1-x_0\rangle=-\gamma_1 \alpha^2\rho_1^2$. Let $L_p$, $H_p$ satisfy \eqref{ineq:bar-Lip-p}. Applying \eqref{ineq:general-smooth} yields 
\[
\begin{aligned}
\chi_1
&=\langle g_0,x_1-x_0\rangle
+\langle g_1-g_0,x_1-x_0\rangle\le-\gamma_1\alpha^2\rho_1^2
+\|g_1-g_0\|_q\|x_1-x_0\|_p\\
&\le-\gamma_1\alpha^2\rho_1^2
+M_1\|x_1-x_0\|_p^2=(M_1-\gamma_1)\alpha^2\rho_1^2\leq M_1\alpha^2\rho_1^2, \quad \text{and}\\
\|s_1\|_q&=\|g_1\|_q
\le\|g_0\|_q+\|g_1-g_0\|_q
\le \alpha\rho_1(\gamma_1+M_1).
\end{aligned}
\]
Using these estimates in \eqref{def:betak+1}, we obtain
\[
\begin{aligned}
\gamma_2& =\max\left\{\theta\gamma_1,\,
\frac{\chi_1+\sqrt{\chi_1^2+2\rho_1^2\|s_1\|_q^2}}
{2\rho_1^2}\right\} \le \max\left\{\theta\gamma_1,\,
\frac{\max\{\chi_1,0\}}{\rho_1^2}
+\frac{\|s_1\|_q}{\sqrt2\rho_1}\right\}\\
& \le
\max\left\{\theta\gamma_1,\,
\alpha^2M_1
+\frac{\alpha(\gamma_1+M_1)}{\sqrt2}\right\}\leq \alpha^2M_1
+\frac{\alpha(\gamma_1+M_1)}{\sqrt2},
\end{aligned}
\]
where the first inequality follows from $\sqrt{a^2+b^2}\le|a|+|b|$, the last follows from $\alpha>\sqrt2$, and $\theta \in (0,1)$. Since the second entry is also at least $\gamma_1$, we obtain $\Gamma_0\leq \alpha^2M_1
+{\alpha(\gamma_1+M_1)}/{\sqrt2}$.
\end{proof}

The next lemma provides upper bounds for $\{\gamma_k\}_{k\geq 2}$ generated by Algorithm \ref{alg:pODA}.

\begin{lemma}
\label{lem:gamma-upper-bound}
Suppose that Assumption~\ref{asp:cvx} holds. Let $\{\gamma_k\}$ be generated by Algorithm \ref{alg:pODA}, and let $r_q$, $\{\ell_k\}$,$\{M_k\}$, $\Gamma_{\mathcal{B}}$, and $\Gamma_0$ be defined in \eqref{def:rq}, \eqref{def:ellk}, \eqref{def:Mk}, \eqref{def:GammaB}, and \eqref{def:gamma0} respectively. Define  $\Gamma_1 := \max\{r_a, \Gamma_{\mathcal{B}}\}$. Then, for all $k\geq 1$,  
\begin{equation} \label{def:gk}
\gamma_{k+1}\le G_k,\qquad G_k:=
\begin{cases}
\Gamma_0, & k=1,\\
\Gamma_0+\Gamma_1
\|(h_2M_2,\ldots,h_kM_k)\|_{r_q},
& k\ge2.
\end{cases}  
\end{equation}
\end{lemma}

\begin{proof}
We prove $\gamma_{k+1}\le G_k$ for all $k\ge 1$ by induction, starting from the established bound
$\gamma_2\le\Gamma_0$.
Fix $k\ge2$ and suppose that $\gamma_k\le G_{k-1}$. Consider the following two cases.
\begin{itemize}
    \item[(i)] Case $k\in\mathcal A$. By invoking the bound $\gamma_{k+1}\le \eta_{k+1}$ established in Lemma \ref{lem:descent-iterations}, together with the definition of $\eta_{k+1}$, the estimate $\ell_k\le h_kM_k$ in \eqref{ineq:upbd-local-estimate}, the inequality $r_a\le \Gamma_1$ implied by the definition of $\Gamma_1$, and the induction hypothesis $\gamma_k\le G_{k-1}$, it follows that
\begin{align*}
    \gamma_{k+1}
    &\le
    \eta_{k+1}
    =\|(\gamma_k,r_a\ell_k)\|_{r_q}\le\bigl\|(
            \Gamma_0 + \Gamma_1 \|(h_2M_2,\ldots,h_{k-1}M_{k-1})\|_{r_q},
            r_a\ell_k)\bigr\|_{r_q} \\
    &\le\bigl\|(
            \Gamma_0 + \Gamma_1 \|(h_2M_2,\ldots,h_{k-1}M_{k-1})\|_{r_q},
            r_ah_k M_k)\bigr\|_{r_q}\le
    \Gamma_0 + \Gamma_1 \|(h_2M_2,\ldots, h_kM_k)\|_{r_q}=G_k.
\end{align*}
% where we adopt the convention that the norm of an empty sequence to be zero.
\item[(ii)] Case $k\in\mathcal B$. By invoking the bound
$\gamma_{k+1}\le\Gamma_{\mathcal B}h_kM_k$ established in Lemma~\ref{lem:gamma-bound-B}, together with the inequality $\Gamma_{\mathcal B}\le\Gamma_1$
implied by the definition of $\Gamma_1$,
we obtain 
\[
\gamma_{k+1} \le \Gamma_{\mathcal B}h_kM_k\le \Gamma_0 + \Gamma_1 \|(h_2M_2,\ldots,h_kM_k)\|_{r_q}=G_k.\] 
\end{itemize}
Combining these two cases completes the proof.
\end{proof}

% The next theorem establish the convergence rate of Algorithm \ref{alg:pODA}.
%

Next, we provide our main theorem which establishes the convergence rate of Algorithm \ref{alg:pODA}.
\begin{theorem}
\label{thm:convex-convergence-rate}
Suppose that Assumption~\ref{asp:cvx} holds. Let
$\{x_k\}$ be generated by Algorithm~1,
% with $q\in[1,2]$, $m\ge2$, $\theta\in(0,1)$, \textcolor{blue}{$\gamma_1>0$, and $c>2(1+e+e/\theta)$}, let $p=q/(q-1)$, 
let $r_q$, $\{\ell_k\}$, $\bar{\rho}$, $\Gamma_{\mathcal B}$, $\Gamma_0$, and $\{G_k\}$ be defined in \eqref{def:rq}, \eqref{def:ellk}, \eqref{eq:bar-rho}, \eqref{def:GammaB}, \eqref{def:gamma0} and \eqref{def:gk}, respectively, and let $\Gamma_1$ be defined in Lemma \ref{lem:gamma-upper-bound}. Then, for any integer $T\ge 1$,
\begin{align}\label{ineq:convex-adaptive-rate}
    f(x_T)-f_*
    \le
    \frac{B_p G_T}{A_T}, \qquad B_p = \|x_0-x_*\|_p^2/2 + {\bar\rho^2}.
\end{align}
Moreover, for $T\ge2$, we further have
\begin{align} \label{eq:main-result}
f(x_T)-f_*
&\le
\frac{2m^m\Gamma_0B_p}{T^m}
+\frac{2^m m^{2m}\Gamma_1L_pB_p}
{T^{3-2/q}}+
\left(
2^{m(1+\nu)}m^{2m}
\Gamma_1H_pCB_p
\right)^{\frac1{1-\nu}}
r_{m,\nu,q}(T),
\end{align}
where $C:=C(m,\nu,q)>0$ is independent of $T$ and is given by
\begin{equation}
\begin{aligned} \label{def:C}
C&=
\begin{cases}
1, & q=2\ \text{or}\ m(1-\nu)=3-\frac{2}{q},\\
\left[
\max\left\{
1,\frac{2/q-1}{|m(1-\nu)-3+2/q|}
\right\}
\right]^{2/q-1}, &
\text{otherwise},
\end{cases}\\
\text{and }\; r_{m,\nu,q}(T) &= \begin{cases}
        T^{-m}, & m < \frac{3-2/q}{1-\nu},\\
        T^{-\frac{3-2/q}{1-\nu}}(\log T)^{\frac{2/q-1}{1-\nu}}, & m=\frac{3-2/q}{1-\nu},\\
        T^{-\frac{3-2/q}{1-\nu}}, & m>\frac{3-2/q}{1-\nu}.
    \end{cases} 
\end{aligned}
\end{equation}
\end{theorem}

\begin{proof}
For $T \ge1$, by the definitions of $\mathcal P_T$, $s_T$, and $c_T$ in \eqref{def:Lyapunov} and \eqref{def:sk-ck-chik}, we have
\begin{align*}
    \mathcal P_T =
    A_Tf_T-c_T+\frac{\|s_T\|_q^2}{2\gamma_{T+1}} =
    A_T(f_T-f_*)
    +\sum_{k=1}^{T}\lambda_k D_f(x_*,x_k)
    +\langle s_T,x_*-x_0\rangle
    +\frac{\|s_T\|_q^2}{2\gamma_{T+1}}.
\end{align*}
It then follows from the convexity of $f$, H\"older's and Young's inequalities that
\begin{align*}
    \mathcal P_T
    \ge
    A_T(f_T-f_*)
    -\|s_T\|_q\|x_0-x_*\|_p
    +\frac{\|s_T\|_q^2}{2\gamma_{T+1}} \ge
    A_T(f_T-f_*)
    -\frac{\gamma_{T+1}}{2}\|x_0-x_*\|_p^2.
\end{align*}
Combining this with
$\mathcal P_T\le {\rho_T^2}\gamma_{T+1}$ in \eqref{eq:potential-budget-bound}, $\gamma_{k+1}\leq G_k$ for all $k\geq 1$ in \eqref{def:gk},
we obtain
\begin{align*}
    A_T(f_T-f_*)
    \le
    \mathcal P_T+\frac{\gamma_{T+1}}{2}\|x_0-x_*\|_p^2 \overset{\eqref{eq:potential-budget-bound}}{\leq}
    \frac{\gamma_{T+1}}{2}
    (\|x_0-x_*\|_p^2+2{\rho_T^2})
    {\le B_p\gamma_{T+1}}
    \overset{\eqref{def:gk}}{\le} B_pG_T,
\end{align*}
where we use $\rho_T\le\bar\rho$ and the definition of $B_p$.
Dividing both sides by $A_T$ yields
\eqref{ineq:convex-adaptive-rate}.
We now assume $T\ge2$ and estimate $G_T$.
For simplicity, we denote
\begin{equation} \label{eq:sim-def}
V_T:=\|(h_k)_{k=2}^T\|_{r_q}, \qquad U_T:=
\|(h_kA_{k-1}^{-\nu})_{k=2}^T\|_{r_q}.   
\end{equation}
First, for $2\leq k\leq T$, applying \eqref{ineq:convex-adaptive-rate} with $T$ replaced by $k-1$ and $k$, respectively, and using $G_{k-1}\le G_k\le G_T$, we have
$\Delta(x_{k-1}) \le B_p G_T/A_{k-1}$ and $\Delta(x_{k}) \le B_p G_T/A_{k}$.
Moreover, by the definition of $M_k$ and the monotonicity of $\{A_k\}$, we have 
    \begin{align*}
        M_k = L_p + H_p\max\{\Delta(x_{k-1}), \Delta(x_k)\}^\nu \le L_p + H_pB_p^\nu G_T^\nu A_{k-1}^{-\nu}.
    \end{align*}
This together with the definition of $G_T$ gives
 \begin{align}
G_T
&=\Gamma_0+\Gamma_1
\|(h_2M_2,\ldots,h_TM_T)\|_{r_q}\le \Gamma_0+\Gamma_1L_p
\|(h_2,\ldots,h_T)\|_{r_q}+\Gamma_1H_pB_p^\nu G_T^\nu
\|(h_2A_1^{-\nu},\ldots,
h_TA_{T-1}^{-\nu})\|_{r_q} \nonumber\\
&\overset{\eqref{eq:sim-def}}{=}
\Gamma_0+\Gamma_1L_pV_T
+\Gamma_1H_pB_p^\nu U_TG_T^\nu\le 2\Gamma_0+2\Gamma_1L_pV_T
+\left(2\Gamma_1H_pB_p^\nu U_T
\right)^{\frac{1}{1-\nu}}, \label{eq:gt-est}
\end{align}
where the first inequality uses Minkowski's inequality, and the last inequality follows from Lemma~\ref{lem:sublinear-gronwall-type} with $a=G_T$, $b=\Gamma_0+\Gamma_1L_pV_T$, and $c=\Gamma_1H_pB_p^\nu U_T$.
We then estimate $V_T$ and $U_T$. Using Lemma~\ref{lem:weights-estimates}, we have 
$A_k\ge c_Ak^m$ and $h_k\le C_hk^{m-2}$, where $c_A=(m(m-1)^{m-1})^{-1}$ and $C_h=2^{m-1}/c_A$ are the same constants as those defined in Lemma~\ref{lem:weights-estimates}.
It follows that
\begin{align} 
V_T
&=\|(h_k)_{k=2}^T\|_{r_q}
\le C_h\|(k^{m-2})_{k=2}^T\|_{r_q}\le C_hT^{m-2+1/r_q}
=C_hT^{m-3+2/q}, \label{eq:vt-bound}\\
U_T
& =\|(h_kA_{k-1}^{-\nu})_{k=2}^T\|_{r_q}\le C_h\left(\frac{2^m}{c_A}\right)^\nu
\|(k^{m(1-\nu)\!-\!2})_{k=2}^T\|_{r_q}, \label{eq:ut-bound}
\end{align}
where we use $1/r_q=2/q-1$, with
$1/r_q=0$ when $q=2$. 
If $r_q < \infty$ (equivalently, $q\in[1,2)$), define $\beta= (m(1-\nu)-2)r_q$. The standard power-sum estimate yields 
\begin{align}\label{ineq:power-sum}
    \sum_{k=2}^T k^\beta \le \begin{cases}
        C_\beta, & \beta < -1,\\
        C_\beta \log T, & \beta = -1,\\
        C_\beta T^{\beta+1}, & \beta > -1,
    \end{cases}
    \qquad
    C_\beta=
\begin{cases}
-{1}/{(\beta+1)}, & \beta<-1,\\
1, & \beta=-1,\\
\max\left\{1,{1}/({\beta+1})\right\},
& \beta>-1.
\end{cases}
\end{align}
Taking the $r_q$-th root of \eqref{ineq:power-sum}
and using \eqref{eq:ut-bound}, we obtain
\begin{align} \label{eq:bound-ut-refined-1}
U_T
\le C_hC_\beta^{1/r_q} \left(\frac{2^m}{c_A}\right)^\nu
\tilde{r}_{m,\nu,q}(T), \qquad \tilde{r}_{m,\nu,q}(T) =\begin{cases}
    1, & m < \frac{3-2/q}{1-\nu}, \\
    (\log T)^{2/q-1}, & m = \frac{3-2/q}{1-\nu}, \\
    T^{m(1-\nu)-3+2/q}& m > \frac{3-2/q}{1-\nu}.
\end{cases}
\end{align}
If $r_q=\infty$ (equivalently, $q=2$), we instead take the maximum of $k^{m(1-\nu) -2}$ in \eqref{eq:ut-bound}.
Consequently,
\begin{align}\label{eq:bound-ut-refined-2}
U_T
\le C_h\left(\frac{2^m}{c_A}\right)^\nu
\tilde{r}_{m,\nu,2}(T), \qquad \tilde{r}_{m,\nu,2}(T)=\begin{cases}
    1, & m \le \frac{2}{1-\nu},\\
    T^{m(1-\nu)-2}, & m > \frac{2}{1-\nu}.
\end{cases}
\end{align}
Combining the two cases in \eqref{eq:bound-ut-refined-1} and \eqref{eq:bound-ut-refined-2}, for each fixed $q\in[1,2]$ with $C$ defined in \eqref{def:C}, we have $U_T \le C_hC \left(2^m/c_A\right)^\nu
\tilde{r}_{m,\nu,q}(T)$. Substituting the bounds on $V_T$ and $U_T$
into the estimate of $G_T$, we obtain
\begin{align*}
G_T
&\overset{\eqref{eq:gt-est}}{\leq} 2\Gamma_0+2\Gamma_1L_pV_T
+\left(2\Gamma_1H_pB_p^\nu U_T
\right)^{\frac{1}{1-\nu}}\\
&\overset{\eqref{eq:vt-bound},\eqref{eq:bound-ut-refined-1},\eqref{eq:bound-ut-refined-2}}{\le} 2\Gamma_0
+2\Gamma_1L_pC_hT^{m-3+2/q}+
\left(2\Gamma_1H_pB_p^\nu C_hC \left(2^m/c_A\right)^\nu
\right)^{\frac{1}{1-\nu}}
\,\tilde{r}_{m,\nu,q}^{\frac{1}{1-\nu}}(T).
\end{align*}
Finally, using \eqref{ineq:convex-adaptive-rate}
and $A_T\ge c_AT^m$, we obtain
\begin{align*}
f(x_T)-f_*
&\le \frac{B_pG_T}{A_T}\le \frac{B_p}{c_A}\bigg[
2\Gamma_0T^{-m}
+2\Gamma_1L_pC_hT^{-(3-2/q)}+
\left(2\Gamma_1H_pB_p^\nu C_hC \left(2^m/c_A\right)^\nu
\right)^{\frac{1}{1-\nu}}
\,T^{-m}\,\tilde{r}_{m,\nu,q}^{\frac{1}{1-\nu}}(T)
\bigg].
\end{align*}
Since $r_{m,\nu,q}(T)
=T^{-m}\tilde r_{m,\nu,q}(T)^{1/(1-\nu)}$,
rearranging the constants and using $m-1\le m$ yields
\eqref{eq:main-result}, which completes the proof.
\end{proof}

\begin{remark}[Upper bounds on $\Gamma_0B_p$]
\label{re:gamma0bp}
We now derive upper bounds on $\Gamma_0B_p$ for the cases with and without calibration of $\gamma_1$. The calibration procedure is described in Calibration \ref{alg:gamma1-calibration}.

\begin{itemize}[leftmargin=1.2em]
\item[(i)] \underline{Without calibration of $\gamma_1$}.
By convexity, $\|g_0\|_q=\alpha\gamma_1\rho_1$, and the definitions of $\Gamma_0$ in \eqref{def:gamma0} and $B_p$ in \eqref{ineq:convex-adaptive-rate}, we have
$\Delta(x_0)\le\alpha\gamma_1B_p\le\alpha\Gamma_0B_p$.
Moreover, \eqref{ineq:convex-adaptive-rate} with $T=1$ gives
$\Delta(x_1)\le\Gamma_0B_p$, since $A_1=1$ and $G_1=\Gamma_0$.
Since $\alpha>1$, the quantity $M_1$ defined in \eqref{def:Mk} satisfies
\[
M_1=L_p+H_p\max\{\Delta(x_0),\Delta(x_1)\}^{\nu}
\le L_p+H_p(\alpha\Gamma_0B_p)^\nu.
\]
The bound $\Gamma_0\le\alpha^2M_1+\alpha(\gamma_1+M_1)/\sqrt2$ in Lemma~\ref{lem:initial-scale-bound} yields
\[
\begin{aligned}
\Gamma_0B_p
& \le
(\tfrac{\alpha}{\sqrt2}\gamma_1+(\alpha^2+\tfrac{\alpha}{\sqrt2}) L_p)B_p
+\alpha^\nu (\alpha^2+\tfrac{\alpha}{\sqrt2}) H_pB_p(\Gamma_0B_p)^\nu \\
&\leq 
2\big(\tfrac{\alpha}{\sqrt2}\gamma_1+(\alpha^2+\tfrac{\alpha}{\sqrt2}) L_p\big)B_p
+\left(2\alpha^\nu (\alpha^2+\tfrac{\alpha}{\sqrt2}) H_pB_p\right)^{\frac1{1-\nu}}.
\end{aligned}
\]
where the second inequality uses Lemma~\ref{lem:sublinear-gronwall-type}.
For arbitrary $\gamma_1>0$, combining this estimate with Remark~\ref{re:barrho} gives 
\begin{align}
&\quad\;\Gamma_0B_p 
\le
2\kappa_0
\left(\tfrac{\alpha}{\sqrt2}\gamma_1+(\alpha^2+\tfrac{\alpha}{\sqrt2}) L_p\right)\|x_0-x_*\|_p^2+
\left(2\alpha^\nu (\alpha^2+\tfrac{\alpha}{\sqrt2}) H_p\kappa_0\|x_0-x_*\|_p^2\right)^{\frac1{1-\nu}} \nonumber\\
&%\!\sim \!
\!=\!\mathcal{O}\!\Bigg(\!
(\gamma_1\!+\!L_p)
\left[
1\!+\!\frac{(L_p\!+\!H_p\Delta(x_0)^\nu)^2}{\gamma_1^2}
\right]\|x_0\!-x_*\|_p^2+\!
\left\{
H_p\!\!
\left[
1\!+\!\frac{(L_p\!+\!H_p\Delta(x_0)^\nu)^2}{\gamma_1^2}
\right]\!\!\|x_0\!-x_*\|_p^2
\right\}^{\!\frac1{1-\nu}}\!
\Bigg) \nonumber\\
&\text{with}\; \kappa_0 = 
\tfrac12+
\max\left\{
\tfrac{(L_p+H_p\Delta(x_0)^\nu)^2}{\alpha^2\gamma_1^2},
\tfrac{C_\theta^2}{\alpha^2-C_\theta^2}
\right\}. \label{eq:gamma0bp-without-calibration}
\end{align}
This inequality follows from the estimate given in Remark~\ref{re:barrho}
with 
\[
\begin{aligned}
B_p
&=
\frac12\|x_0-x_*\|_p^2
+
\max\left\{
\frac{\|g_0\|_q^2}{\alpha^2\gamma_1^2},
\frac{C_\theta^2\|x_0-x_*\|_p^2}{\alpha^2-C_\theta^2}
\right\}\\
&\le
\left[
\frac12+
\max\left\{
\frac{(L_p+H_p\Delta(x_0)^\nu)^2}{\alpha^2\gamma_1^2},
\frac{C_\theta^2}{\alpha^2-C_\theta^2}
\right\}
\right]\|x_0-x_*\|_p^2=\kappa_0\|x_0-x_*\|_p^2.
\end{aligned}
\]
\item[(ii)] \underline{With calibration of $\gamma_1$}.
If $\gamma_1$ is selected by Calibration \ref{alg:gamma1-calibration}, 
then
$\|g_0\|_q^2/(2\Delta(x_0))\le\gamma_1<2(L_p+H_p\Delta(x_0)^\nu)$ in
\eqref{pro:gamma1}
and $\Delta(x_1)\le\Delta(x_0)$ hold. Hence, the quantity $M_1$ defined in \eqref{def:Mk} satisfies $M_1=L_p+H_p\Delta(x_0)^\nu$, and applying Lemma~\ref{lem:initial-scale-bound} together with $\gamma_1<2M_1$ yields $\Gamma_0\le(\alpha^2+3\alpha/\sqrt2)M_1$.
Furthermore, the convexity of $f$ and \eqref{ineq:general-smooth} imply $\Delta(x_0)\le M_1\|x_0-x_*\|_p^2$. Substituting this into  $M_1$ yields
\[
\begin{aligned}
M_1\|x_0-x_*\|_p^2
& \le
L_p\|x_0-x_*\|_p^2
+H_p\|x_0-x_*\|_p^2
\left(M_1\|x_0-x_*\|_p^2\right)^\nu \\
&\leq
2L_p\|x_0-x_*\|_p^2+
\left(2H_p\|x_0-x_*\|_p^2\right)^{\frac1{1-\nu}},
\end{aligned}
\]
where the second inequality uses Lemma~\ref{lem:sublinear-gronwall-type}.
We then obtain
\begin{equation}
\label{eq:gamma0bp-with-calibration}
\begin{aligned}
\Gamma_0B_p
&\le 
\left(\alpha^2+\frac{3\alpha}{\sqrt2}\right)M_1B_p=\left(\alpha^2+\frac{3\alpha}{\sqrt2}\right)(L_p+H_p\Delta(x_0)^\nu)B_p\\
&\leq 
\left(\alpha^2+\frac{3\alpha}{\sqrt2}\right)
\left[
\frac12+
\max\left\{
\frac4{\alpha^2},
\frac{C_\theta^2}{\alpha^2-C_\theta^2}
\right\}
\right]
\left[
2L_p\|x_0-x_*\|_p^2+
\left(2H_p\|x_0-x_*\|_p^2\right)^{\frac1{1-\nu}}
\right]\\
& %\sim 
=\mathcal{O}\left(
L_p\|x_0-x_*\|_p^2
+
\left(H_p\|x_0-x_*\|_p^2\right)^{\frac1{1-\nu}}
\right),
\end{aligned}
\end{equation}
where the last inequality follows from the estimate given in Remark~\ref{re:barrho}
with 
\[
B_p
\le
\left[
\frac12+
\max\left\{
\frac4{\alpha^2},
\frac{C_\theta^2}{\alpha^2-C_\theta^2}
\right\}
\right]\|x_0-x_*\|_p^2.
\]
\end{itemize}
% we have
% \[
% \Gamma_0B_p
% \lesssim_{\alpha,\theta,\nu}
% (\gamma_1+L_p)\|x_0-x_*\|_p^2
% +\left(H_p\|x_0-x_*\|_p^2\right)^{\frac1{1-\nu}}.
% \]
% \end{itemize}
% % Thus, calibration removes the additional dependence on
% % $(L_p+H_p\Delta(x_0)^\nu)/\gamma_1$ from this upper bound, while $\gamma_1$ still appears linearly in the first term.
Estimate \eqref{eq:gamma0bp-without-calibration} holds for arbitrary $\gamma_1>0$ and retains explicit dependence on $\gamma_1$ and $\Delta(x_0)$. Calibration yields \eqref{eq:gamma0bp-with-calibration}, which removes this dependence and the additional factors involving $(L_p+H_p\Delta(x_0)^\nu)^2/\gamma_1^2$, thereby improving the dependence on $L_p$ and $H_p$ in the upper bound on $\Gamma_0B_p$.
\end{remark}

\begin{remark}[Convergence rate under Calibration \ref{alg:gamma1-calibration}]
\label{re:calibrated-convergence-rate}
This remark aims to confirm the convergence rate stated earlier in \eqref{re:clear-rate}.
Under the assumptions of Theorem~\ref{thm:convex-convergence-rate}, with $0\le\nu<1$, choose $\gamma_1$ using Calibration~\ref{alg:gamma1-calibration}. Remarks~\ref{re:barrho} and~\ref{re:gamma0bp} then bound $B_p$ and $\Gamma_0B_p$, respectively. Substituting these estimates into \eqref{eq:main-result} yields the following explicit convergence bound for every integer $T\ge2$:
\[
\begin{aligned}
f(x_T)-f_*
&\le
\frac{2m^m}{T^m}
\left(\alpha^2+\frac{3\alpha}{\sqrt2}\right)
\left[
\frac12+\max\left\{
\frac4{\alpha^2},
\frac{C_\theta^2}{\alpha^2-C_\theta^2}
\right\}
\right]
\left[
2L_p\|x_0-x_*\|_p^2+
\left(2H_p\|x_0-x_*\|_p^2\right)^{\frac1{1-\nu}}
\right]
\\
&\quad+
\frac{2^m m^{2m}\Gamma_1}{T^{3-2/q}}
\left[
\frac12+\max\left\{
\frac4{\alpha^2},
\frac{C_\theta^2}{\alpha^2-C_\theta^2}
\right\}
\right]
L_p\|x_0-x_*\|_p^2
\\
&\quad+
\left\{
2^{m(1+\nu)}m^{2m}\Gamma_1C
\left[
\frac12+\max\left\{
\frac4{\alpha^2},
\frac{C_\theta^2}{\alpha^2-C_\theta^2}
\right\}
\right]
\right\}^{\frac1{1-\nu}}
\left(H_p\|x_0-x_*\|_p^2\right)^{\frac1{1-\nu}}
r_{m,\nu,q}(T),
\end{aligned}
\]
where $C$ and $r_{m,\nu,q}$ are defined in \eqref{def:C}. In addition, recall that $\Gamma_1=\max\{r_a,(4+a_q^{-1})(1+\alpha)^2\}$, defined in Lemma~\ref{lem:gamma-upper-bound}, is independent of $\gamma_1$.
This calibration sharpens the upper bounds on both $\Gamma_0B_p$ and $B_p$, and consequently improves the dependence of the convergence bound on $L_p$ and $H_p$.
\end{remark}

\bibliographystyle{abbrv}
\bibliography{ref}

\appendix

\section*{Appendix}

\section{Proof of the $q$-norm expansion inequality}

\begin{proof}[Proof of Lemma~\ref{lem:expan-vq}]
When $u=0$, \eqref{ineq:V-expan} holds clearly. Suppose that
$\|u\|_q>0$. Lemma \ref{lem:expan-omegaq} gives
\begin{align}\label{ineq:u+vq-expansion}
    \|u+v\|_q^q
    \le
    \|u\|_q^q + q\sum_{i=1}^n \omega_q(u_i)v_i + 2^{2-q}\|v\|_q^q = \|u\|_q^q + \frac{q}{\|u\|_q^{2-q}}\langle J_q(u), v\rangle
    +2^{2-q}\|v\|_q^q,
\end{align}
where the last equality uses the definitions of $\omega_q$ and $J_q$. Recall from \eqref{ppty:J} that $\|J_q(u)\|_p=\|u\|_q$. By this and H\"older's inequality, one has $|\langle J_q(u),v\rangle|\le\|J_q(u)\|_p\|v\|_q=\|u\|_q\|v\|_q$. We now derive a nonnegative upper bound for the first two terms on the right-hand side of the above inequality:
\begin{align*}
    \|u\|_q^q
    +\frac{q}{\|u\|_q^{2-q}}\langle J_q(u),v\rangle
    &\le
    \left|
        \|u\|_q+\frac{\langle J_q(u),v\rangle}{\|u\|_q}
    \right|^q=
    \left(
        \|u\|_q^2
        +2\langle J_q(u),v\rangle
        +\frac{\langle J_q(u),v\rangle^2}{\|u\|_q^2}
    \right)^{q/2}\\
    &\le
    \left(
        \|u\|_q^2
        +2\langle J_q(u),v\rangle
        +\|v\|_q^2
    \right)^{q/2},
\end{align*}
where the first inequality is due to the convexity of $t\mapsto |t|^q$. Combining this with \eqref{ineq:u+vq-expansion} yields
\begin{align}\label{ineq:power-q-expansion-intermediate}
    \|u+v\|_q^q
    \le(
        \|u\|_q^2
        +2\langle J_q(u),v\rangle
        +\|v\|_q^2)^{q/2}
    +2^{2-q}\|v\|_q^q.
\end{align}
By H\"older's inequality with conjugate exponents
$2/q$ and $2/(2-q)$, one has that $a_1b_1+a_2b_2\le (a_1^{2/q}+a_2^{2/q})^{q/2}(b_1^{2/(2-q)}+b_2^{2/(2-q)})^{(2-q)/2}$ for $a_1, b_1, a_2, b_2\geq 0$. Fix any $c\in(0,1)$. By \eqref{ineq:power-q-expansion-intermediate} and applying this H\"older's inequality with
\begin{align*}
    a_1 = (
        \|u\|_q^2
        +2\langle J_q(u),v\rangle
        +\|v\|_q^2)^{q/2},\quad a_2 = \frac{2^{2-q}c^{q/2}\|v\|_q^q}
     {(1-c^{r_q})^{(2-q)/2}},\quad b_1 = c^{q/2},\quad b_2 = (1-c^{r_q})^{(2-q)/2},
\end{align*}
we obtain
\begin{align*}
c^{q/2}\|u+v\|_q^q
&\le
c^{q/2}
(
\|u\|_q^2
+2\langle J_q(u),v\rangle
+\|v\|_q^2
)^{q/2}
+(1-c^{r_q})^{(2-q)/2}
\frac{2^{2-q}c^{q/2}\|v\|_q^q}
     {(1-c^{r_q})^{(2-q)/2}}
\\
&\le
\bigg(
\|u\|_q^2
+2\langle J_q(u),v\rangle
+\|v\|_q^2
+\bigg(
\frac{2^{2-q}c^{q/2}\|v\|_q^q}
     {(1-c^{r_q})^{(2-q)/2}}
\bigg)^{2/{q}}
\bigg)^{q/2}
\\
&=
\bigg(
\|u\|_q^2
+2\langle J_q(u),v\rangle
+\bigg(
1+\frac{2^{2/{r_q}}c}{(1-c^{r_q})^{1/{r_q}}}
\bigg)\|v\|_q^2
\bigg)^{q/2},
\end{align*}
where we use the definition ${r_q}=q/(2-q)$.
Raising both sides to the power $2/q$ and dividing by $2$ proves \eqref{ineq:V-expan} as desired.
\end{proof}

\section{Generalized local cocoercivity}\label{apx:cc}

\begin{lemma}
\label{lem:generalized-local-cocoercivity}
Suppose that $f:\mathbb R^n\to\mathbb R$ is convex and continuously
differentiable and attains its minimum value $f_*$. Define
$\Delta(x):=f(x)-f_*$, and fix an arbitrary optimal solution $x_*$.
Let $p\in[2,\infty]$ and $q\in[1,2]$ satisfy $1/p+1/q=1$, and let $x_c\in\mathbb R^n$ be arbitrary. Fix any $R>0$ satisfying $R\ge\|x_*-x_c\|_p$, and define
\[
    \mathcal C_p
    :=\{x\in\mathbb R^n:\|x-x_c\|_p\le R\},
    \qquad
    \bar{\mathcal C}_p
    :=\{x\in\mathbb R^n:\|x-x_c\|_p\le3R\}.
\]
Then, for any $L_p,H_p\ge0$ and $\nu\in[0,1)$, the Lipschitz-type condition
\begin{align}\label{eq:generalized-lip-local}
    \|\nabla f(u)-\nabla f(v)\|_q
    \le
    \left(
        L_p+H_p\max\{\Delta(u),\Delta(v)\}^{\nu}
    \right)\|u-v\|_p
    \qquad
    \forall u,v\in\bar{\mathcal C}_p
\end{align}
implies the cocoercivity-type condition
\begin{align}\label{eq:generalized-local-cocoercivity}
    \|\nabla f(x)-\nabla f(y)\|_q^2
    \le
    2\left(
        L_p+2^\nu H_p\max\{\Delta(x),\Delta(y)\}^{\nu}
    \right)D_f(x,y)
    \qquad
    \forall x,y\in\mathcal C_p.
\end{align}
\end{lemma}

\begin{proof}
Suppose that \eqref{eq:generalized-lip-local} holds. We now prove
\eqref{eq:generalized-local-cocoercivity}. Fix any $x,y\in\mathcal C_p$. For convenience, denote
\begin{align*}
&g_x:=\nabla f(x),\quad g_y:=\nabla f(y),\quad \Phi_{x,y}:=\max\{\Delta(x),\Delta(y)\},\\
&M_{x,y}:= L_p + H_p\,\Phi_{x,y}^\nu,\quad \widetilde{M}_{x,y}:= L_p + 2^\nu H_p\,\Phi_{x,y}^\nu.
\end{align*}
The conclusion is immediate if $g_x= g_y$. Henceforth, suppose that $g_x\neq g_y$, which implies $\Phi_{x,y}>0$. For any $z\in\mathcal C_p$ and any $v\in\mathbb R^n$ with
$\|v\|_p=1$, we define $\varphi_{z,v}(s):=\Delta(z-sv)$ for all $s\ge0$. Taking the derivative, we have
$\varphi_{z,v}'(s)=-\langle\nabla f(z-sv),v\rangle$. By this,
Hölder's inequality, and \eqref{eq:generalized-lip-local}, one has
\begin{align}\label{eq:generalized-phi-est}
    |\varphi_{z,v}'(s)-\varphi_{z,v}'(r)|
    \le \|\nabla f(z-sv)-\nabla f(z-rv)\|_q  \le
    (L_p+H_p\max\{\varphi_{z,v}(s),\varphi_{z,v}(r)\}^{\nu})|s-r|
\end{align}
for all $s,r\ge0$ such that $z-sv,z-rv\in\bar{\mathcal{C}}_p$. This one-dimensional Lipschitz-type condition will be used repeatedly below.

We first bound $\|g_y\|_q$ in terms of the objective gap. Claim that
\begin{align}\label{eq:generalized-gradf-bdd}
    \|g_y\|_q^2 \le 2M_{x,y}\Phi_{x,y}.
\end{align}
This claim holds clearly if $\|g_y\|_q =0$. Otherwise, let $w\in\mathbb R^n$ be such that $\|w\|_p=1$ and
$\langle g_y,w\rangle=\|g_y\|_q>0$. Since $\nabla f(x_*)=0$ and $x_*\in\mathcal C_p$, it follows from
\eqref{eq:generalized-lip-local} that
\[
    \|g_y\|_q \le M_{x,y}\|y-x_*\|_p \le 2M_{x,y}R.
\]
Thus, setting $s_{x,y}:=\|g_y\|_q/M_{x,y}$, we have
$\{y-sw:0\le s\le s_{x,y}\}\subset\bar{\mathcal C}_p$, since
\[
    \|y-sw-x_c\|_p
    \le \|y-x_c\|_p+s\|w\|_p
    \le R+s_{x,y}
    \le 3R
    \qquad \forall s\in[0,s_{x,y}].
\] Then, applying \eqref{eq:generalized-phi-est} with $(z,v)=(y,w)$ and $r=0$ and using $\varphi_{y,w}(0)=\Delta(y)\le\Phi_{x,y}$, we obtain
\begin{align}
\varphi_{y,w}'(s) %& %\le \varphi_{y,w}'(0) + (L_p+H_p\max\{\varphi_{y,w}(s),\varphi_{y,w}(0)\}^{\nu})s \nonumber\\
    & \le \varphi_{y,w}'(0) + (L_p+H_p\max\{\varphi_{y,w}(s),\Phi_{x,y}\}^{\nu})s\qquad\forall s\in[0,s_{x,y}].\label{ineq:upper-bound-varphi-yw1}
\end{align}
We next prove $\varphi_{y,w}(s)\le\Phi_{x,y}$ for all $s\in[0,s_{x,y}]$ and relate the coefficient of $s$ in \eqref{ineq:upper-bound-varphi-yw1} to $M_{x,y}$. Notice from the definition of $\varphi_{y,w}$ and the choice of $w$ that $\varphi_{y,w}'(0) = -\|g_y\|_q<0$. We show by contradiction that
$\varphi_{y,w}(s)<\varphi_{y,w}(0)$ for all $s\in(0,s_{x,y}]$. Otherwise, there exists an
$s_0\in(0,s_{x,y}]$ such that
$\varphi_{y,w}(s_0)=\varphi_{y,w}(0)$ and
$\varphi_{y,w}(s)<\varphi_{y,w}(0)$ for all $s\in(0,s_0)$.
Therefore, $\varphi_{y,w}(s)\le\varphi_{y,w}(0)\le\Phi_{x,y}$ for all $s\in[0,s_0]$. By this and \eqref{ineq:upper-bound-varphi-yw1}, one has
\begin{align*}
\varphi_{y,w}'(s) \le \varphi_{y,w}'(0) + (L_p+H_p\Phi_{x,y}^{\nu})s = \varphi_{y,w}'(0)+M_{x,y}s =-\|g_y\|_q+M_{x,y} s\qquad\forall s\in[0,s_0].%\label{ineq:upper-bound-varphi-yws}
\end{align*}
Integrating this inequality over $s\in[0,s_0]$, we obtain
$\varphi_{y,w}(s_0)\le\varphi_{y,w}(0)-\|g_y\|_q\,s_0
+M_{x,y}s_0^2/2<\varphi_{y,w}(0)$,
where the strict inequality follows from
$s_0\in(0,s_{x,y}]$ and $s_{x,y}=\|g_y\|_q/M_{x,y}$.
This contradicts $\varphi_{y,w}(s_0)=\varphi_{y,w}(0)$.
Hence,
$\varphi_{y,w}(s)\le\varphi_{y,w}(0)\le\Phi_{x,y}$
for all $s\in[0,s_{x,y}]$. Combining this with \eqref{ineq:upper-bound-varphi-yw1}, we obtain
\[
    \varphi_{y,w}'(s) \le \varphi_{y,w}'(0) + (L_p+H_p\Phi_{x,y}^{\nu})s =-\|g_y\|_q+M_{x,y}s
    \qquad
    \forall s\in[0,s_{x,y}].
\]
Integrating this inequality over $s\in[0,s_{x,y}]$, we obtain 
\begin{align*}%\label{eq:generalized-delta-est}
\varphi_{y,w}(s_{x,y})
    \le \varphi_{y,w}(0)-\|g_y\|_q\,s_{x,y}+\frac{M_{x,y}}{2}s_{x,y}^2 = \varphi_{y,w}(0) - \frac{\|g_y\|_q^2}{2M_{x,y}} \le \Phi_{x,y} - \frac{\|g_y\|_q^2}{2M_{x,y}}.
\end{align*}
Together with $\varphi_{y,w}(s_{x,y})\ge0$, which follows from
the definition of $\varphi_{y,w}$, this yields the claim
\eqref{eq:generalized-gradf-bdd}.

We are now ready to prove \eqref{eq:generalized-local-cocoercivity}. Let $u\in\mathbb R^n$ be such that $\|u\|_p=1$ and
$\langle g_x-g_y,u\rangle=\|g_x-g_y\|_q$. By $x,y\in\mathcal{C}_p\subset\bar{\mathcal C}_p$, one has from \eqref{eq:generalized-lip-local} that
\begin{align}\label{eq:generalized-direct-bdd}
    \|g_x-g_y\|_q \le M_{x,y} \|x-y\|_p \le 2M_{x,y}R.
\end{align}
Since $g_x\neq g_y$, it then follows that
$M_{x,y}>0$, and hence $\widetilde M_{x,y}\ge M_{x,y}>0$.
Set $t_{x,y}:=\|g_x-g_y\|_q/\widetilde M_{x,y}$. By \eqref{eq:generalized-direct-bdd}, we have $t_{x,y}\le2R$.
Thus,
$\{x-su:0\le s\le t_{x,y}\}\subset\bar{\mathcal C}_p$, since
\[
    \|x-su-x_c\|_p
    \le \|x-x_c\|_p+s\|u\|_p
    \le R+t_{x,y}
    \le3R
    \qquad \forall s\in[0,t_{x,y}].
\]
Then, applying \eqref{eq:generalized-phi-est} with $(z,v)=(x,u)$
and $r=0$, and using $\varphi_{x,u}(0)=\Delta(x)\le\Phi_{x,y}$,
we obtain
\begin{align}\label{ineq:upper-bound-varphi-xu}
    \varphi_{x,u}'(s)
    &\le \varphi_{x,u}'(0)
    +(L_p+H_p\max\{\varphi_{x,u}(s),\Phi_{x,y}\}^{\nu})s
    \qquad \forall s\in[0,t_{x,y}].
\end{align}
Claim that
\begin{align}\label{eq:generalized-phix-est}
    \varphi_{x,u}(s)\le2\Phi_{x,y}
    \qquad \forall s\in[0,t_{x,y}].
\end{align}
Suppose for contradiction that this bound fails.
Since $\varphi_{x,u}(0)\le\Phi_{x,y}<2\Phi_{x,y}$,
continuity implies that there exists an
$s_0\in(0,t_{x,y})$ such that
$\varphi_{x,u}(s_0)=2\Phi_{x,y}$ and $\varphi_{x,u}(s)\le2\Phi_{x,y}$ for all $s\in[0,s_0]$.
Combining this with \eqref{ineq:upper-bound-varphi-xu} yields $\varphi_{x,u}'(s)\le\varphi_{x,u}'(0)+\widetilde M_{x,y}s$ for all $s\in[0,s_0]$. By the choice of $u$ and Hölder's inequality, one has
\[
    \varphi_{x,u}'(0)
    =-\langle g_x,u\rangle
    =-\|g_x-g_y\|_q-\langle g_y,u\rangle
    \le\|g_y\|_q-\|g_x-g_y\|_q.
\]
For convenience, denote
\[
U(t) := \Phi_{x,y} + (\|g_y\|_q-\|g_x-g_y\|_q)\,t + \frac{\widetilde M_{x,y}}{2}\,t^2\qquad\forall t\in[0,t_{x,y}].
\]
Then, by $\varphi_{x,u}(s)\le2\Phi_{x,y}$ for all $s\in[0,s_0]$
and integrating \eqref{ineq:upper-bound-varphi-xu} over $s\in[0,s_0]$,
we obtain
\begin{align}\label{ineq:upbd-varphi-xu-s0}
    \varphi_{x,u}(s_0)\le\varphi_{x,u}(0) + \varphi_{x,u}'(0)\,s_0 + (L_p+2^\nu H_p\Phi_{x,y}^{\nu})\,\frac{s_0^2}{2} \le U(s_0).
\end{align}
On the other hand, by the definitions of $U$ and $t_{x,y}$, we have
\begin{align*}
    U(0) = \Phi_{x,y},\qquad U(t_{x,y})
    =\Phi_{x,y}
    +\frac{\|g_y\|_q\,\|g_x-g_y\|_q}{\widetilde M_{x,y}}
    -\frac{\|g_x-g_y\|_q^2}{2\widetilde M_{x,y}}\le\Phi_{x,y}
    +\frac{\|g_y\|_q^2}{2\widetilde M_{x,y}}
    \le2\Phi_{x,y},
\end{align*}
where the first inequality follows by completing the square,
and the second inequality follows from \eqref{eq:generalized-gradf-bdd}
and $M_{x,y}\le\widetilde M_{x,y}$. By convexity of $U$ and the fact that $s_0\in[0,t_{x,y})$, one has
\[
    U(s_0)
    \le\left(1-\frac{s_0}{t_{x,y}}\right)U(0)
    +\frac{s_0}{t_{x,y}}U(t_{x,y})
    \le\left(1+\frac{s_0}{t_{x,y}}\right)\Phi_{x,y}
    <2\Phi_{x,y},
\]
which along with \eqref{ineq:upbd-varphi-xu-s0} implies $\varphi_{x,u}(s_0)\le U(s_0)<2\Phi_{x,y}$. This contradicts $\varphi_{x,u}(s_0)=2\Phi_{x,y}$. Hence, the claim \eqref{eq:generalized-phix-est} holds.

Combining \eqref{eq:generalized-phix-est} with
\eqref{ineq:upper-bound-varphi-xu}, we obtain
\begin{align}\label{ineq:upbd-varphi-xu-txy-}
    \varphi_{x,u}'(s)
    \le \varphi_{x,u}'(0)
    +\left(L_p+2^\nu H_p\Phi_{x,y}^{\nu}\right)s
    =\varphi_{x,u}'(0)+\widetilde M_{x,y}s
    \qquad \forall s\in[0,t_{x,y}].    
\end{align}
Define $h(z):=f(z)-\langle g_y,z\rangle$ for all $z\in\mathbb R^n$.
Since $h$ is convex and $\nabla h(y)=0$, $y$ is a global
minimizer of $h$. Moreover, by the definitions of $h$ and $\varphi_{x,u}$, one has
\begin{align}\label{eq:generalized-h-prop}
    h(x)-h(y)=D_f(x,y),
    \quad
    \nabla h(x)=g_x-g_y,\quad h(x-su)-h(x)
    =\varphi_{x,u}(s)-\varphi_{x,u}(0)
    +s\langle g_y,u\rangle.
\end{align}
Integrating \eqref{ineq:upbd-varphi-xu-txy-} over $s\in[0,t_{x,y}]$
and using $\varphi_{x,u}'(0)=-\langle g_x,u\rangle$, we obtain
\begin{align*}
    h(x-t_{x,y}u)-h(x)
    \le -t_{x,y}\langle g_x-g_y,u\rangle
    +\frac{\widetilde M_{x,y}}{2}t_{x,y}^2 =-\|g_x-g_y\|_q\,t_{x,y}
    +\frac{\widetilde M_{x,y}}{2}t_{x,y}^2.
\end{align*}
Since $y$ is a global minimizer of $h$, we have
$h(y)\le h(x-t_{x,y}u)$. Therefore,
\begin{align*}
    D_f(x,y)
    =h(x)-h(y) \ge \|g_x-g_y\|_q\,t_{x,y}
    -\frac{\widetilde M_{x,y}}{2}t_{x,y}^2
    =\frac{\|g_x-g_y\|_q^2}{2\widetilde M_{x,y}}.
\end{align*}
Multiplying both sides by $2\widetilde M_{x,y}$ proves
\eqref{eq:generalized-local-cocoercivity} as desired.
\end{proof}

\section{Proof of the calibration bound}
\label{app:calibration-bound}

\begin{proof}[Proof of \eqref{pro:gamma1}]
Since $\gamma_1$ is accepted and $f(x_1)\ge f_*$,
the sufficient-decrease condition gives $\frac{\|g_0\|_q^2}{2\gamma_1}
\le \Delta(x_0)$,
which proves the lower bound.

For the upper bound, we first show that the rejected trial
$x_1(\gamma_1/2)$ belongs to $\mathcal S_p$ defined in \eqref{def:bounded-set}.
By convexity, for any minimizer $x_*$, we have $\Delta(x_0)
\le \langle g_0,x_0-x_*\rangle
\le \|g_0\|_q\|x_0-x_*\|_p$.
Combining \eqref{ppty:J} with the lower bound above yields
\[
\begin{aligned}
\|x_1(\gamma_1/2)-x_0\|_p
=\frac{2\|g_0\|_q}{\gamma_1}\le \frac{4\Delta(x_0)}{\|g_0\|_q}\le 4\|x_0-x_*\|_p.
\end{aligned}
\]
On the other hand, \eqref{eq:bar-rho} and
$\alpha>C_\theta>4$ imply
\[
\alpha\bar\rho
\ge
\frac{\alpha C_\theta}
{\sqrt{\alpha^2-C_\theta^2}}
\|x_0-x_*\|_p
>4\|x_0-x_*\|_p.
\]
Hence $x_1(\gamma_1/2)\in\mathcal S_p$.
Second, we show that every trial
$\gamma\ge\max\{\gamma_1/2,L_p+H_p\Delta(x_0)^\nu\}$
satisfies the sufficient-decrease condition.
Define
\[
z(t)=x_0+\frac{t}{\gamma}J_q(-g_0),
\qquad t\in[0,1].
\]
Since $\mathcal S_p$ is convex and contains
$x_0$ and $x_1(\gamma_1/2)$,
we have $z(t)\in\mathcal S_p$ for all $t\in[0,1]$.
Suppose that $f(z(t))\ge f(x_0)$
for some $t\in(0,1]$.
Since
\[
\left.\frac{\mathrm d}{\mathrm d t}f(z(t))\right|_{t=0}
=-\frac{\|g_0\|_q^2}{\gamma}<0,
\]
there exists a first $t_0\in(0,1]$ such that
$f(z(t_0))=f(x_0)$ and
$f(z(t))<f(x_0)$ for all $t\in(0,t_0)$.
Thus, for $t\in[0,t_0]$,
$\Delta(z(t))\le\Delta(x_0)$.
By \eqref{ineq:general-smooth}, we have
\[
\|\nabla f(z(t))-g_0\|_q
\le (L_p+H_p\Delta(x_0)^\nu)\|z(t)-x_0\|_p
=\frac{(L_p+H_p\Delta(x_0)^\nu)t}{\gamma}\|g_0\|_q.
\]
Integrating along the segment and using
\eqref{ppty:J}, we obtain
\[
\begin{aligned}
f(z(t_0))-f(x_0)
\le
-\frac{t_0}{\gamma}\|g_0\|_q^2
+\frac{(L_p+H_p\Delta(x_0)^\nu)t_0^2}{2\gamma^2}\|g_0\|_q^2\le
-\frac{t_0}{2\gamma}\|g_0\|_q^2<0,
\end{aligned}
\]
where the last inequality uses $\gamma\ge (L_p+H_p\Delta(x_0)^\nu)$
and $t_0\le1$. This is a contradiction.
Consequently, $f(z(t))\le f(x_0)$
throughout $[0,1]$. Applying the same
estimate at $t=1$ gives
\[
\begin{aligned}
f(x_1(\gamma))
\le f(x_0)
-\frac{\|g_0\|_q^2}{\gamma}
+\frac{(L_p+H_p\Delta(x_0)^\nu)\|g_0\|_q^2}{2\gamma^2}\le f(x_0)
-\frac{\|g_0\|_q^2}{2\gamma}.
\end{aligned}
\]
Thus, every $\gamma\ge\max\{\gamma_1/2,L_p+H_p\Delta(x_0)^\nu\}$
is accepted.
By the calibration procedure, the trial
$\gamma_1/2$ is rejected. Therefore, we obtain $\gamma_1<2(L_p+H_p\Delta(x_0)^\nu)$.
This completes the proof.
\end{proof}

\section{Technical lemmas}

\begin{lemma}\label{lem:weights-estimates}
    Let $c_A = (m(m-1)^{m-1})^{-1}$, and $C_h=2^{m-1}/c_A$. For $m \ge 2$ and $k\ge 2$, the weights $\lambda_k$ and $A_k$ in Algorithm~\ref{alg:pODA} satisfy
    \begin{align} \label{eq:results-3in}
        1\le \frac{\lambda_k}{\lambda_{k-1}} \le e, \qquad c_A k^m\le A_k \le k^m, \qquad \frac{k^{m-2}}{2(m-1)^{2m-2}}\le h_k := \frac{\lambda_k^2}{2A_{k-1}}\le C_h k^{m-2}.
    \end{align}
\end{lemma}
\begin{proof}
    We prove the results given in \eqref{eq:results-3in} one by one. First, for $k\ge 2$, using the update of $\lambda_k$ in \eqref{def:Ak-lambdak}, the first inequality can be verified directly by
    \begin{align*}
        1\leq \frac{\lambda_k}{\lambda_{k-1}} = \left(\frac{k+m-2}{k+m-3}\right)^{m-1} 
        % = \left(1 + \frac{1}{k+m-3}\right)^{m-1}
        \leq \left(1 + \frac{1}{m-1}\right)^{m-1}\leq e.
    \end{align*}
    Second, for every $k\ge 1$, we have 
    $
        \frac{k}{m-1} \le 1 + \frac{k-1}{m-1}=\lambda_k^{1/(m-1)}\le k.
    $
    Hence, taking the $m-1$ power, summing $\lambda_i$ from $i=1$ to $k$, and using the update form of $A_k$ in \eqref{def:Ak-lambdak}, we obtain
    \begin{align*}
        \frac{k^{m-1}}{(m-1)^{m-1}} \le \lambda_k \le k^{m-1},\qquad 
\frac{1}{(m-1)^{m-1}}\sum_{i=1}^k i^{m-1}
\le A_k
\le\sum_{i=1}^k i^{m-1}.
    \end{align*}
    Using the elementary inequality $k^m/m \leq \sum_{i=1}^{k} i^{m-1}\leq k^m$ for $m\geq 1$, we directly have $c_A k^m \le A_k \le k^m$ by using the definition $c_A = (m(m-1)^{m-1})^{-1}$.
    % Multiplying the lower bound by $(m-1)^{m-1}/m$ and summing $\lambda_i$ from $i=1$ to $k$ proves $c_A k^m \le A_k \le k^m$. 
    %
    Finally, for $k\ge 2$, the inequality $k-1\ge k/2$ gives
    \begin{align*}
        h_k = \frac{\lambda_k^2}{2A_{k-1}} \le \frac{k^{2m-2}}{2c_A(k-1)^m} \le \frac{2^{m-1}}{c_A} k^{m-2} = C_h k^{m-2}.
    \end{align*}
    For the reverse order estimate, using $\lambda_k \ge (k/(m-1))^{m-1}$ and $A_{k-1} \le (k-1)^m$, we obtain
    \begin{align*}
        h_k = \frac{\lambda_k^2}{2A_{k-1}} \ge \frac{k^{2m-2}}{2(m-1)^{2m-2}(k-1)^m} \ge \frac{k^{m-2}}{2(m-1)^{2m-2}}.
    \end{align*}
    This completes the proof.    
\end{proof}

\begin{lemma}\label{lem:expan-omegaq}
Let $q\in[1,2)$, and define $\omega_q(t):=\operatorname{sgn}_0(t)|t|^{q-1}$ for all $t\in\R$.
% \footnote{We adopt the convention $0^0=1$.} 
Then,
\begin{align}\label{ineq:power-q-expansion}
    \|u+v\|_q^q
    \le
    \|u\|_q^q
    +q\sum_{i=1}^n \omega_q(u_i)v_i
    +2^{2-q}\|v\|_q^q
\end{align}
holds for all $u,v\in\R^n$.
\end{lemma}

\begin{proof}
Define $\omega_q(t):=\operatorname{sgn}_0(t)|t|^{q-1}$ for all $t\in\R$. We first show that, for all $q\in[1,2)$,
\begin{align}\label{ineq:diff-jq}
    |\omega_q(a)-\omega_q(b)|
    \le 2^{2-q}|a-b|^{q-1}
    \qquad\forall a,b\in\R.
\end{align}
To prove the claim, we consider the following three cases.
\begin{itemize}[leftmargin=1.8em]
    \item[(i)] If $ab>0$, the subadditivity of $t\mapsto t^{q-1}$ on $\mathbb{R}_+$ gives $|\omega_q(a)-\omega_q(b)|
=\big||a|^{q-1}-|b|^{q-1}\big|
\le \big||a|-|b|\big|^{q-1}
=|a-b|^{q-1}$, which further implies \eqref{ineq:diff-jq} using $q\in[1,2)$;
    % where the inequality is due to $1-t^{q-1}\le(1-t)^{q-1}$ for all $t\in[0,1]$.
    \item[(ii)] If $ab<0$, the concavity of $t\mapsto t^{q-1}$ on $\mathbb{R}_+$ gives $|\omega_q(a)-\omega_q(b)|
    =|a|^{q-1}+|b|^{q-1}
    \le
    2^{2-q}(|a|+|b|)^{q-1}
    =
    2^{2-q}|a-b|^{q-1}$;
    \item[(iii)] If $ab=0$, then $|\omega_q(a)-\omega_q(b)|\leq|a-b|^{q-1}$, which further implies \eqref{ineq:diff-jq}.
\end{itemize}

Combining all these cases establishes \eqref{ineq:diff-jq}. Moreover, for any $q\in[1,2)$, the function $t\mapsto |t|^q$
is locally absolutely continuous, and its derivative equals
$q\omega_q(t)$ almost everywhere.
Thus, for all $a,b\in\R$,
\begin{align*}
    |a+b|^q-|a|^q-q\omega_q(a)b
    &=
    q\int_0^1
    (\omega_q(a+tb)-\omega_q(a))b\,\mathrm{d}t\le
    q\int_0^1
    |\omega_q(a+tb)-\omega_q(a)|\,|b|\,\mathrm{d}t\\
    &\overset{\eqref{ineq:diff-jq}}{\le}
    q\,2^{2-q}|b|^q
    \int_0^1 t^{q-1}\,\mathrm{d}t=2^{2-q}|b|^q.
\end{align*}
Applying it with $a=u_i$ and $b=v_i$ and summing over
$i=1,\ldots,n$ yields the desired result \eqref{ineq:power-q-expansion}.
\end{proof}

\begin{lemma}\label{lem:sublinear-gronwall-type}
    Let $a, b, c \ge 0$, $\nu \in [0,1)$. If
    $
        a \le b + ca^\nu,
    $
    then
    $
        a \le 2b + (2c)^{\frac{1}{1-\nu}}.
    $
\end{lemma}
\begin{proof}
    By the weighted Young's inequality, we have 
    \begin{align*}
        c a^{\nu} = \frac12 \bigl[(2c)^{1/(1-\nu)}\bigr]^{1-\nu}a^\nu \le \frac12 \Big((1-\nu)(2c)^{\frac{1}{1-\nu}}+\nu a\Big) \le \frac12 \Big((2c)^{\frac{1}{1-\nu}} + a\Big).
    \end{align*}
    Substituting this into the inequality $a \le b + ca^\nu$ and rearranging gives the claim.
\end{proof}

\end{document}

%% file: math_commands.tex
\usepackage{amsmath,amsfonts,bm}

\def\1{\bm{1}}

\DeclareMathAlphabet{\mathsfit}{\encodingdefault}{\sfdefault}{m}{sl}
\SetMathAlphabet{\mathsfit}{bold}{\encodingdefault}{\sfdefault}{bx}{n}

\newcommand{\R}{\mathbb{R}}

\DeclareMathOperator*{\argmin}{arg\,min}